\documentclass[12pt,a4paper,twoside]{article}
\usepackage{amsmath}
\usepackage{amsfonts}
\usepackage{amssymb}
\usepackage{amsthm}
\usepackage{enumerate}
\usepackage{amscd}
\usepackage{amsrefs}
\usepackage{hyperref}
\usepackage{verbatim}
\usepackage{calligra}
\usepackage{mathrsfs}
\usepackage{xcolor}
\usepackage[all,cmtip]{xy}
\usepackage[upright]{fourier}
\usepackage[version=3]{mhchem}

\fontfamily{Liberation Serif}

\hypersetup{citebordercolor={0 1 0}}
\newcommand{\et}{\text{e}}
\newcommand{\Ht}{\text{H}}

\newcommand{\Kt}{\text{K}}
\newcommand{\Et}{\text{E}}

\newcommand{\ts}{\text{Spec}}

\newcommand{\hit}{\text{ht}}
\newcommand{\grd}{\text{grade}}
\newcommand{\As}{\text{Ass}}

\newcommand{\as}{\text{ass}}
\newcommand{\dep}{\text{Depth}}
\newcommand{\dime}{\text{dim}}

\newcommand{\map}{\mathfrak{p}}
\newcommand{\mam}{\mathfrak{m}}
\newcommand{\maq}{\mathfrak{q}}
\newcommand{\maa}{\mathfrak{a}}
\newcommand{\mab}{\mathfrak{b}}

\newcommand{\mad}{\mathfrak{d}}
\newcommand{\man}{\mathfrak{n}}

\newcommand{\mn}{\mathbb{N}}

\newcommand{\ext}{\text{Ext}}
\newcommand{\homm}{\text{Hom}}

\newcommand{\dsum}{\bigoplus}

\newcommand{\ten}{\bigotimes}
\newcommand{\tenr}{\bigotimes_R}

\newcommand{\ins}{\bigcap}
\newcommand{\uni}{\bigcup}
\newcommand{\unil}{\bigcup\limits}

\newcommand{\Frac}{\text{Frac}}
\newcommand{\Ra}{R(a^{1/2})}

\DeclareMathAlphabet{\mathcalligra}{T1}{calligra}{g}{f}

\begin{document}
   
    \date{}
    \title{\normalsize \textbf{Reduction of some homological conjectures to  excellent unique factorization domains}}
    \markboth{\small \textsc{$\ $Reduction of some homological conjectures  to ...}}{\small  \textsc{Ehsan Tavanfar}}         
      \author{\normalsize \textsc{Ehsan Tavanfar\footnote{This research was in part supported by a grant from IPM and supported in part by NSF grant DMS 1162585.}$\ $}\\ \\ \normalsize \textsc{In the memory of my father, Manouchehr Tavanfar,}\\ \normalsize \textsc{ who passed away at the time of writing  this paper}\\}
    \maketitle          
    \begin{abstract}
       \textsc{Abstract}. 
         In this article, applying the quasi-Gorenstein analogous of  the Ulrich's deformation of certain Gorenstein rings  we show that some homological conjectures, including the Monomial Conjecture,  Big Cohen-Macaulay Algebra Conjecture as well as the Small Cohen-Macaulay Conjecture reduce to the  excellent unique factorization domains. Some other reductions of the Monomial Conjecture are also proved. We, moreover, show that certain almost complete intersections adhere the Monomial Conjecture.
    \end{abstract}     
    
    {\let\thefootnote\relax}
    {\let\thefootnote\relax}
    
    \newtheorem{thm}{Theorem}[section]
    \theoremstyle{Definition}    
    \newtheorem{defi}[thm]{Definition}
    \theoremstyle{Definition and Remark}
    \newtheorem{defi-rem}[thm]{Definition and Remark}
    \newtheorem{defis-rem}[thm]{Definitions and Remark}
    \newtheorem{defis-rems}[thm]{Definitions and Remarks}
    \newtheorem{defi-Nots}[thm]{Definition and Notations}
    \newtheorem{defi-Not}[thm]{Definition and Notation}
    \newtheorem{notas}[thm]{Notations}
    \theoremstyle{Lemma}
    \newtheorem{lem}[thm]{Lemma}
    \theoremstyle{remark}
    \newtheorem{rem}[thm]{Remark}
    \theoremstyle{Corollary}
    \newtheorem{cor}[thm]{Corollary}
    \newtheorem{exam}[thm]{Example}
    \newtheorem{counterexam}[thm]{Counterexample}
    \newtheorem{prop}[thm]{Proposition}
    \newtheorem{ques}[thm]{Question}
    \newtheorem{discuss}[thm]{Discussion}
    \newtheorem{conj}[thm]{Conjecture}
    \normalsize  
    
    \section{Introduction}

    Throughout this article, $(R,\mathfrak{m})$ denotes a commutative Noetherian local ring of dimension $d$ with identity where $\mathfrak{m}$  denotes the unique maximal ideal of $R$.  
    
    In the early seventies, in \cite{HochsterContracted}, Hochster proposed the following two equivalent conjectures and proved them for the equicharacteristic rings. 
      
       \textbf{Monomial Conjecture:}  Suppose that $x_1,\ldots,x_d$ is a system of parameters for $R$. Then, $x_1^t\cdots x_d^t\notin (x_1^{t+1},\ldots,x_d^{t+1})$ for each $t\ge 1$.
    
    \textbf{Direct Summand Conjecture:} Suppose that $A$ is a regular local ring and $R$ is a module finite extension of $A$. Then the inclusion $A\rightarrow R$ splits, in the category of $A$-modules.
       
       The foregoing conjectures have been around for decades in the mixed characteristic case. In the sequel we succinctly review certain substantial advances in the study of the aforementioned conjectures. In the eighties, in his extraordinary paper \cite{HochsterCanonical}, among other things, Hochster introduced the following further equivalent forms of the Monomial Conjecture.
     
    \textbf{Canonical Element Conjecture:} Suppose that $(F_\bullet,\partial_\bullet)$ is a free resolution of $R/\mam$ and consider the natural projection $\pi:F_d\rightarrow \text{syz}^d_R(R/\mam)$ which defines an element, $$\epsilon_R:=[\pi]\in \ext^d_{R}\big(R/\mam,\text{syz}^d_R(R/\mam)\big).$$
    Set, $\eta_R:=\phi(\epsilon_R),$ wherein, $$\phi:\ext^d_{R}\big(R/\mam,\text{syz}^d_R(R/\mam)\big)\rightarrow \lim\limits_{\overset{\longrightarrow}{n\in \mn}}\ext^d_{R}\big(R/\mam^n,\text{syz}^d_R(R/\mam)\big)=\Ht^d_\mam\big(\text{syz}^d_R(R/\mam)\big),$$ is the natural map. Then $\eta_R\neq 0$.
    
    \textbf{CE Property, also known as the Canonical Element Conjecture:} Assume that $F_\bullet$ is minimal resolution of $R/\mam$. Consider the Koszul complex $\Kt_\bullet(\mathbf{x},R)$ of a system of parameters $\mathbf{x}$ of $R$. Then for any lifting, $\phi:\Kt_\bullet(\mathbf{x},R)\rightarrow F_\bullet$, of the natural surjection $R/\mathbf{(x)}\rightarrow R/\mam$, we have $\phi_d\neq 0$.

    In \cite{HochsterCanonical}, Hochster, moreover shows that the mentioned equivalent forms of the Monomial Conjecture imply the following conjecture:

    \textbf{Improved New Intersection Conjecture:} Suppose that $F_\bullet:0\rightarrow F_s\rightarrow F_{s-1}\rightarrow \cdots \rightarrow F_0\rightarrow 0$ is a complex of finite free modules such that $\ell\big(\Ht_i(F_\bullet)\big)<\infty$ for each $i>1$ and $\Ht_0(F_\bullet)$ has a minimal generator whose generated $R$-module has finite length. Then $s\ge d$.
    
    The important point here is that a slightly  stronger version of the Improved New Intersection Conjecture, wherein $\Ht_0(F_0)$ has finite length as well, has been settled in full generality by Roberts  (see,\cite{RobertsLeTheorem}). The validity of this stronger form, which is called the New Intersection Theorem, itself established a bunch of other old conjectures, most notably the Bass's question regarding the Cohen-Macaulayness of rings having a non-zero finite module of finite injective dimension (see, \cite[Chapter 9.]{BrunsHerzogCohenMacaulay} and \cite{PeskineSzpiroDimension}). Originally the Improved New Intersect Conjecture is, implicitly, stated in \cite{EvansGriffithTheSyzygy} among the proof of \cite[Theorem 1.1]{EvansGriffithTheSyzygy}. Subsequently, miraculously, in \cite{DuttaOnTheCanonical} Dutta showed that the improved new intersection is, in fact, an another equivalent form of the Monomial Conjecture. Recently, in \cite{DuttaTheMonomialII } Dutta  also proved that all of the equivalent forms of the Monomial Conjecture are in addition equivalent to the following conjecture which has its origin, again, in \cite{EvansGriffithTheSyzygy} (see, \cite[Proposition 1.6.]{EvansGriffithTheSyzygy}).
    
    \textbf{Order Ideal Conjecture:} Assume that $R$ is a regular ring and $M$ is an $R$-module. Suppose that $x$ is a minimal generator of the $i$th-syzygy of $M$ for some $i>0$. Then the ideal, $$\mathcal{O}_{\text{syz}^i(M)}(x)=\{f(x):f\in \homm_R\big(\text{syz}^i(M),R\big)\},$$ has height greater than or equal to $i$.
    
      An $R$-module $M$ (algebra $A$)  is said to be a balanced big Cohen-Macaulay module ( big Cohen-Macaulay algebra) if $\mam M\neq M$ ($\mam A\neq A$) and any system of parameters $\mathbf{x}$ of $R$ is a regular sequence on $M$ (on $A$).  The existence of balanced big Cohen-Macaulay module (big Cohen-Macaulay algebra), which is a long standing conjecture, settles the Monomial Conjecture affirmatively. In \cite{HochsterTheEquicharacteristic}  Hochster showed that any equicharacteristic ring has a balanced big Cohen-Macaualy module.  Moreover, in \cite{HochsterHunekeInfinite }  Hochster and Huneke proved that if $R$ is an excellent domain of prime characteristic, then,  $R^+$, the integral closure of 	$R$ in the algebraic closure of its fraction field is a big Cohen-Macaulay algebra. Unfortunately, using the trace map, we can deduce that the analogous statement about $R^+$ in equal characteristic zero does not hold provided $\dim(R)\ge 3$. In spite of this obstruction, in \cite{HochsterHunekeApplications } Hochster and Huneke established the existence of big Cohen-Macaulay algebras in the equicharacteristic case. 
      
      Although the Monomial Conjecture is quite easy in dimension less than or equal to two,  but not in dimension three. In \cite{HeitmannDirect} Heitmann provided a very complicated proof of the validity of Monomial Conjecture in dimension three of mixed characteristic.  In fact he proved more, that is, if $R$ is a mixed characteristic excellent normal domain of dimension three then $\Ht^2_\mam(R^+)$ is almost zero, i.e. $R^+$ is an almost Cohen-Macaulay $R$-algebra\footnote{Here, almost zero means that each element of $\Ht^2_\mam(R^+)$ is killed by elements of arbitrary small valuation, with respect to a fixed rank one valuation which is non-negative on $R$ and positive on $\mam$.}.  Quite obviously, the existence of almost Cohen-Macaulay algebra  is stronger than the existence of big Cohen-Macaulay algebras. In general  the Monomial Conjecture follows from the existence of almost Cohen-Macaulay algebra(see, e.g. \cite[Proposition 1.3]{RobertsSinghAnnihilators}). In \cite{HochsterBigCohenMacaulay } when $R$ has dimension three and mixed characteristic Hochster constructed a big Cohen-Macaulay algebra using the Heitmann's result regarding the almost Cohen-Macaulayness of $R^+$. Recently, in \cite{BhattacharyyaExistence}, Bhattacharyya  showed that the Big Cohen-Macaulay Algebra Conjecture follows from the existence of almost big Cohen-Macaulay algebras.  The question whether $R^+$ is always almost Cohen-Macaulay in mixed characteristic is open in dimension greater than or equal to $4$. Although the Monomial Conjecture, has been proved in equal characteristic zero but the question whether $R^+$ is almost Cohen-Macaulay in equal characteristic zero is of self interest and is open in dimension greater than or equal to $3$. But the three dimensional equal characteristic zero case for the Segre product of Cohen-Macaulay $\mathbb{N}$-graded rings has been investigated  in  \cite{RobertsSinghAnnihilators }, where the authors apply the Albanese varieties to deduce that the image of non-top local cohomologies of $R$  are almost zero in $R^+$. Recently, in \cite{RobertsLocalCohomology } Roberts showed that the image of non-top local cohomologies of $R$ in $R^+$ are almost zero provided $R$ is a  Segre product of $\mathbb{N}$-graded Cohen-Macaulay rings of mixed characteristic. Moreover, recently, applying the Almost Purity Theorem, in \cite{ShimomotoAnApplication} Shimomoto proved that $R$ has a big Cohen-Macaulay algebra provided $R$  has mixed characteristic $p>0$ such that $R[1/p]$ is an \'etale extension of $A[1/p]$ wherein $A$ is a Noether normalization of $R$. 
      
      A finitely generated big  Cohen-Macaulay $R$-module is said to be a small Cohen-Macaulay module or maximal Cohen-Macaulay module. Hochster's Small Cohen-Macaulay Conjecture states that every complete ring has a small Cohen-Macaulay Module. This conjecture is even open in  three dimensional equicharacteristic rings.
          
           In Section 2  we prove that for the validity of the Monomial Conjecture it suffices to take into account only systems of parameters of $R$ which form a part of a minimal basis of $\mam$. This fact will be used in Section 4 to deduce the main result of that section. Here, it is worth to point out that, since for proving this reduction we descend to a finite free extension of $R$ with complete intersection fibers,  we can merge some of  other known reductions of the Monomial Conjecture with our reduction. For instance, it suffices to confine our attention only on such system of parameters $\mathbf{x}:=x_1,\ldots,x_d$ for almost complete intersections. In view of  \cite{DuttaTheMonomial} the Monomial Conjecture reduces to the almost complete intersections.
      
           In Section 3 we show that  the aforementioned homological conjectures reduce to the class of excellent (homomorphic image of regular) unique factorization domains.  The main  idea of the proof is the fact that the Ulrich's deformation of certain Gorenstein rings to the class of unique factorization domains, developed in \cite{UlrichGorenstein }, has a quasi-Gorenstein counterpart. But prior to applying this reduction, for a normal domain $R$ with canonical ideal $\omega_R$, we endow $R\dsum \omega_R$ with a ring structure such that it is a quasi-Gorenstein domain which is a complete intersection in codimension 1 (recall that the trivial extension is never a domain). For the Monomial Conjecture with the aid of  Hochster's first general grade reduction technique developed in \cite{HochsterProperties}, we deduce that the reduction to excellent unique factorization domains is dimensionwise in the sense that the Monomial Conjecture for  $d$-dimensional local rings descends to  $d$-dimensional excellent unique factorization local  domains. It is noteworthy to mention that in the light of  Heitmann's remarkable paper, \cite{HeitmannCharacterization }, loosely speaking, any complete ring with depth at least 2 is a completion of a unique factorization domain! But, unfortunately, Heitmann's unique factorization domains are  neither excellent nor  catenary \footnote{In fact Heitmann's construction of such unique factorization domains provided us with an example of a non-catenary local unique factorization domain. Prior to this counterexample, it was conjectured that any local unique factorization domain is catenary.}. For instance, his unique factorization domains are not supposed to have a maximal Cohen-Macaulay module or even a canonical module. So the excellence property can be important from this point of view. Besides the known good properties of unique factorization domains, they may have even more tacit important properties as we will point out them in the sequel. The first example of a non-Cohen-Macaulay unique factorization domain was constructed by virtue of the invariant theory, in Bertin's paper \cite{BertinAnneaux }. In fact at that time,  Bertin's counterexample settled negatively  Samuel's conjecture regarding the Cohen-Macaulayness of  unique factorization domains. In spite of the existence of this counterexample, complete unique factorization domains at least in equal characteristic zero, have very good properties.  For instance, as stated in \cite[Page 539]{LipmanUnique } (see, also  \cite{HartshorneOgusOnTheFactoriality }), if $R$ is a complete equicharacteristic zero unique factorization domain of depth $\ge 3$, then $R$ satisfies $S_3$. Moreover, if the residue field is algebraically closed (char 0), then the  depth condition is superfluous. Furthermore, in view of \cite[Page 540]{LipmanUnique }, any complete unique factorization domain with algebraically closed residue field  of characteristic zero is Cohen-Macaulay provided $\dim(R)\le 4$.   So it might be, fairly, more helpful, if we can take more steps and reduce any of the homological conjectures to the class of complete unique factorization domains.  
        
         In Section 4, we investigate  the homological conjectures for  almost complete intersections satisfying $\mam^2\subseteq (\mathbf{x})$ for some system of parameters $\mathbf{x}$ of $R$. In mixed characteristic we furthermore assume that $x_1=p$, the residual characteristic of $R$. We will see that in this case $R$ has a balanced big Cohen-Macaulay algebra.  In the sequel, we explain our motivation behind this investigation.  As we mentioned before, the Monomial Conjecture is proved for equicharacteristic rings of arbitrary dimension and for (at most) $3$-dimensional rings of any characteristic. Due to the Goto's paper \cite{GotoOnTheAssociated }, another subclass of rings which adhere the Monomial Conjecture is Buchsbaum rings.  In fact any system of parameters $\mathbf{x}$ of a Buchsbaum ring satisfies, $\mam \big(\{\mathbf{x}\}^{\lim}_R\big)\subseteq (\mathbf{x})$, wherein $\{\mathbf{x}\}^{\lim}_R=\unil_{t\in \mn}\big((x_1^{t+1},\ldots,x_d^{t+1}):x_1^t\cdots x_d^t\big)$. This inclusion immediately settles the Monomial Conjecture. Now set $\mathcal{A}$ to be the class of rings consisting of rings $R$ satisfying, $$\mam^2\big( \{\mathbf{x}\}^{\lim}_R\big)\subseteq (\mathbf{x}),$$ for some system of parameters $\mathbf{x}$ of $R$. Note that $\mathcal{A}$ contains, properly, the class of Buchsbaum rings (see, \cite[Theorem 5.12(ii)]{TavanfarTousiAStudy}). In order to see why this inclusion is proper, note that in the light of \cite{GotoANote }, there exists a non-Buchsbaum quasi-Buchsbaum ring $R$ such that there are exactly two non-zero non-top local cohomologies of $R$ which both of them are $R/\mam$-vector spaces. Then applying the argument of the proof of \cite[Proposition 1.]{RobertsFontaine } we can conclude that $R$ belongs to the class $\mathcal{A}$. Note also that by \cite[Theorem 5.12.(i)]{TavanfarTousiAStudy}, roughly speaking, $\mathcal{A}$ is a subclass of generalized Cohen-Macaulay rings. Quite obviously, if the Monomial Conjecture is valid for any system of parameters $\mathbf{x}$ of any ring $R$ with $\mam^2\subseteq (\mathbf{x})$, then any member of  $\mathcal{A}$ satisfies the Monomial Conjecture. Although, at this time, we could not prove the validity of the Monomial Conjecture in this general setting, but we succeeded to deduce our desired statement for almost complete intersections satisfying $\mam^2\subseteq (\mathbf{x})$ with the extra assumption that $x_1=p$ in the mixed characteristic $p>0$. For justifying  our investigation in this restricted case, from another point of view, recall that in the light of Dutta's result \cite[1.2 Proposition]{DuttaTheMonomial} in conjunction with \cite[(6.1) Theorem]{HochsterCanonical} Monomial Conjecture reduces to the case where $R$ is an almost complete intersection of mixed characteristic $p>0$ and $x_1,\ldots,x_d$ is a system of parameters for $R$ with $x_1=p$.  The key ingredient of our proof is the fact that those  almost complete intersections have multiplicities up to two in the non-Cohen-Macaulay case. We give an example of a non-Cohen-Macaulay ring $R\in \mathcal{A}$ such that its multiplicity is strictly greater than two (so, $R$ is not an almost complete intersection). 
         
         In view of the Monomial Conjecture, we end the paper with a question regarding existence of a particular non-commutative $R$-algebra.

    \section{\large{Reduction of the Monomial Conjecture to  systems of parameters which are part of a minimal basis of $\mam$}}
    
    This section is devoted to prove that for the validity of the Monomial Conjecture it suffices  only to check  systems of parameters of the ring which can be extended to a minimal set of generators of the maximal ideal. This  result will be used in the  Section 4 to prove that certain almost complete intersections satisfy the Monomial Conjecture.          
    \begin{rem} \label{FirstRemark}
       Let $a\in R$. Then the free $R$-module $R\dsum R$ acquires a ring structure via the following rule,
         $$(r,s)(r^\prime,s^\prime)=(rr^\prime+ss^\prime a,rs^\prime+r^\prime s).$$
 We use the notation $R(a^{1/2})$ to denote the forgoing ring structure of $R\dsum R$. In fact it is easily seen that the map, $R(a^{1/2})\rightarrow R[X]/(X^2-a)$, which takes $(r,s)$ to $(sX+r)+(X^2-a)R[X]$ is an isomorphism of $R$-algebras. We are given   the  extension map $R\rightarrow R(a^{1/2})$ by the rule $r\mapsto (r,0)$ which turns $\Ra$ into a free $R$-module with the basis $\{(1,0),(0,1)\}$. Consequently this  extension is an integral extension of $R$ and it is subject to the following properties which all are easy to verify.  
         \begin{enumerate}
           \item[(i)]  $\dime(R)=\dime\big(\Ra\big)$ and $a$ has a square root in $\Ra$, namely $(0,1)$.
           \item[(ii)] If $a\in \mam$ then $\Ra$ is a local ring with unique maximal ideal $\mam\dsum R$.
           \item[(iii)] \label{RootOfParameterElement} If $a,x_2,\ldots,x_d$ is a system of parameters of $R$ then $a^{1/2},x_2,\ldots,x_d$ is a system of parameters for $\Ra$. 
           \item[(iv)]  
           
           If $\varphi:R\rightarrow S$ is an $R$-algebra wherein, $\varphi(a)$ has a square root in $S$, then there exists a natural induced $R$-algebra homomorphism $\Psi:\Ra\rightarrow S$  defined by $(r,s)\mapsto \varphi(r)+\varphi(s)\varphi(a)^{1/2}$, extending $\varphi$. Moreover, assuming that $R$ and $S$ are domains, $\Psi$ is injective if and only if $\varphi$ is injective and $a$ does not have a square root in $\Frac(R)$.
           
           \item[(v)] If $R$ is a  domain, then $\Ra$ is a domain if and only if $a$ does not have a square root in $\text{Frac}(R)$.  For one implication we note that if $a$ does not have a square root in $K:=\Frac(R)$, then $K[x]/(x^2-a)$ is a field extension of $R$. Therefore, by applying the preceding part we get an embedding $\Ra\rightarrow K[x]/(x^2-a)$ which shows that $\Ra$ is an integral domain. For the reverse implication note that if $a$ has a square root $r/s\in \Frac(R)$ then the identity $(r,s)(r,-s)=0$ shows that $R(a^{1/2})$ is not a domain.           
         \end{enumerate}  
         
       \end{rem}
       
        The  following remark is in fact the proof of Corollary \ref{WeCanAssumeMinimalGenerator}(i).
        
        \begin{rem}  \label{ProofOfCorollary}
       
         \label{SquareOfSequence} Let $x_1,\ldots,x_l$ be a sequence of elements of $R$ contained in the maximal ideal of $R$.  We, inductively, construct the local ring $(R_i,\mathfrak{m}_i)$ by taking a square root of $x_{i}$ in $R_{i-1}$. Then in $R_l$ we have, $$x_i^{1/2}=(\underset{0\text{-th coordinate}}{\underbrace{0}},0,\ldots,0,1,0,\ldots,\underset{(2^l-1)\text{-th coordinate}}{\underbrace{0}}),$$ whose $2^{(i-1)}$-th coordinate is $1$ and others are zero. 
        
        \item[(i)] Let $1\le j\le l$ and $0\le k\le 2^l-1$. We denote the element $(0,\ldots,0,\underset{k-\text{th coordinate}}{\underbrace{1}},0,\ldots,0)$ of $R_l$ by $e_k$. Then we have, $$e_kx_{j}^{1/2}=\begin{cases}e_{k+2^{j-1}}, & (j-1)\text{-\ th\ digit\ of\ }k\text{\ in\ base\ }2\text{ is\ 0} \\ x_je_{k-2^{j-1}}, & (j-1)\text{-\ th\ digit\ of\ }k\text{\ in\ base\ }2\text{ is\ 1}.\end{cases}$$ In order to see why this is the case we induct on the least natural number $s\ge j$ such that $k\le 2^s-1$. In the case where $s=j$ it is easily seen that the  $(j-1)$-th digit of $k$ in its $2$-th base representation is $0$ (is $1$) if and only if $k\le 2^{j-1}-1$ ($k\ge 2^{j-1}$). So an easy use of the multiplication rule of the ring $R_{j}:=R_{j-1}\bigoplus R_{j-1}$ proves the claim (Recall that $R_j$ is subring of $R_l$). Now assume that $s>j$. Then we have, $$e_kx_j^{1/2}=\Big(0,\ldots,\underset{2^{s-1}-1\text{-th\ coordinate}}{\underbrace{0}},\big(0,\ldots,0,\underset{(k-2^{s-1})\text{-th\ coordinate}}{\underbrace{1}},0,\ldots,0\big)x_j^{1/2}\Big).$$	
        
        Now set $k^\prime:=k-2^{s-1}$. Note that the $(j-1)$-th digit of the base $2$ representation of $k$ and $k^\prime$ are equal. Consequently the statement follows from our inductive hypothesis.
        
       \item[(ii)] We are going to show that for each $1\le i\le l$ the projection map $\tau_{2^{(i-1)}}:\mam_l^2\rightarrow R$, which is the projection to the $(2^{(i-1)})$-th coordinate,  is not surjective. In the case where $l=1$ we have $\mam_l^2=(\mam^2+x_1R)\bigoplus \mam$. So, we assume that $l\ge 2$ and the statement is true for smaller values of $l$. Then,
       \begin{equation}
       \label{Combinatorics}
        \mam_l^2=\big(\mam_{l-1}^2+x_lR_{l-1}\big)\bigoplus \mam_{l-1}.       
       \end{equation}
       
        Now, if $i=l$ then $2^{l-1}$-th coordinate  of $\mam_l^2$ is just the first coordinate of, $$\mam_{l-1}=\mam\bigoplus R\bigoplus \cdots\bigoplus R.$$ Hence, clearly, $\tau_{2^{(l-1)}}$ is not surjective. On the other hand if $i\lneq l$ then by our inductive hypothesis $\tau_{2^{(i-1)}}:\mam_{l-1}^2\rightarrow R$ is not surjective which, in the light of the equality (\ref{Combinatorics}), implies the statement immediately.
        
        \item[(iii)] In continuation of our investigation in the previous part, we need to show, also, that the projection map $\tau_{2^{(i-1)}}:x_j^{1/2}R_{l-1}\rightarrow R$ is not surjective unless $i=j$ ($1\le i\le l-1$ and $1\le j\le l-1$). Let $(r_k)_{0\le k\le 2^{l-1}-1}\in R_{l-1}$. Then we have,
        
         \begin{center}
          $(r_k)_{0\le k\le 2^{l-1}-1}x_j^{1/2}= \sum\limits_{k=0}^{2^{l-1}-1}r_ke_{k}x_j^{1/2}=\linebreak\sum\limits_{\substack{k=0\\ (j-1)-\text{th\ digit\ of\ }k\text{\ in\ base\ }2\text{\ is\ }0}}^{2^{l-1}-1}r_ke_{k+2^{j-1}}+\sum\limits_{\substack{k=0\\ (j-1)-\text{th\ digit\ of\ }k\text{\ \ in\ base\ }2\text{\ is\ }1}}^{2^{l-1}-1}r_kx_je_{k-2^{j-1}}$.
         \end{center}
        
        Thus if $i<j$ then evidently $\tau_{2^{i-1}}$ is not surjective. On the other hand if $i>j$ and there exits  $0\le k\le 2^{l-1}-1$ with $k+2^{j-1}=2^{i-1}$ then $k=2^{i-1}-2^{j-1}$ which after a straightforward computation shows that the $(j-1)$-th digit of $k$ in base $2$ is $1$. This proves the non-surjectivity of  $\tau_{2^{i-1}}$.
        
        \item[(iv)]  By means of the arguments of the forgoing part we can, directly, conclude that, $$x^{1/2}_i\notin  \mam^2_{l-1}+(x^{1/2}_1,\ldots,\widehat{x^{1/2}_i},\ldots,x_{l-1}^{1/2},x_l)R_{l-1},\ (i\lneq l)$$ otherwise we must have  $\tau_{2^{(i-1)}}:(x^{1/2}_j)R_{l-1}\rightarrow R$ for some $1\le j\le l-1$ and $j\neq i$ or  $\tau_{2^{(i-1)}}:\mam^2_{l-1}\rightarrow R$ is surjective.
        
      \end{rem}

    The above remark yields the following corollary, as claimed before. In the second part of the subsequent corollary the important point is the equality $\text{embdim}(R)=d+u$ which equals to the embedding dimension of the regular local ring $V[[X_2,\ldots,X_d,Y_1,\ldots,Y_u]]$.
    
    \begin{cor} \label{WeCanAssumeMinimalGenerator}
       \begin{enumerate}
         \item[(i)] The Monomial Conjecture holds if and only if  every system of parameters $x_1,\ldots,x_d$ which is a part of a minimal basis for the maximal ideal satisfies the Monomial Conjecture. 
         \item[(ii)] For the Monomial Conjecture, without loss of generality, we can assume that, $$R=\Big((V,p^{1/2})[[X_2,\ldots,X_d,Y_1,\ldots,Y_u]]\Big)/I,$$ and the system of parameters $x_1,\ldots,x_d$ is the image of $p^{1/2},X_2,\ldots,X_d$ in $R$ where the discrete valuation ring $(V,p^{1/2})$ is a subring  of $R$ and $\text{embdim}(R)=d+u$.
       \end{enumerate}  
        \begin{proof}
          (i) Let $x_1,\ldots,x_d$ be a system of parameters for $R$. Consider the ring $R_d$ as in Remark \ref{ProofOfCorollary}. It suffices to show that $x_1^{1/2},\ldots,x_d^{1/2}$ is a part of a minimal basis for the unique maximal ideal $\mam_d$ of $R_d$. Let $(\alpha_1,\beta_1),\ldots,(\alpha_d,\beta_d)\in R_d=R_{d-1}\bigoplus R_{d-1}$ be such that, $$\sum\limits_{k=1}^{d-1}(\alpha_k,\beta_k)x_k^{1/2}+(\alpha_d,\beta_d)\underset{=x_d^{1/2}}{\underbrace{(0_{R_{d-1}},1_{R_{d-1}})}}\in \mam_d^2=\big(\mam_{d-1}^2+(x_d)R_{d-1}\big)\bigoplus \mam_{d-1}.$$
           Then by a simple computation we get \begin{equation}
             \label{FirstIdentity}
             \sum\limits_{k=1}^{d-1}\alpha_kx_k^{1/2}+\beta_{d}x_d\in \mam_{d-1}^2+x_dR_{d-1},
           \end{equation}
              and, 
            \begin{equation}
            \label{SecondIdentity}
            \sum\limits_{k=1}^{d-1}\beta_kx_k^{1/2}+\alpha_d\in \mam_{d-1}.
            \end{equation} 
            So the identity (\ref{SecondIdentity}) yields $\alpha_d\in \mam_{d-1}$  and thence $(\alpha_d,\beta_d)\in \mam_{d}$. Moreover,  for each $1\le k\le d-1$  we must have $(\alpha_k,\beta_k)\in \mam_d=\mam_{d-1}\dsum R_{d-1}$, otherwise  we get $\alpha_i\notin \mam_{d-1}$ for some $1\le i\le d-1$ which in view of the identity (\ref{FirstIdentity}) yields $x_i^{1/2}\in \mam_{d-1}^{2}+(x_1^{1/2},\ldots,\widehat{x_i^{1/2}},\ldots,x_{d-1}^{1/2},x_d)R_{d-1}$ violating  Remark \ref{ProofOfCorollary} (iv). Consequently, $(\alpha_k,\beta_k)\in \mam_d$ for each $1\le k\le d$. This implies that $x_1^{1/2}+\mam_d^2,\ldots,x_d^{1/2}+\mam_d^2$ is a linearly independent subset of $\mam_d/\mam_d^2$ over $R_d/\mam_d$ and thence $x_1^{1/2},\ldots,x_d^{1/2}$ is part of a minimal basis for $\mam_d$.
             
             (ii) By \cite[(6.1) Theorem]{HochsterCanonical }, for the validity of the Monomial Conjecture it suffices to verify it only for systems of parameters $x_1,\ldots,x_d$ of $R$ wherein $x_1=p$ is the residual characteristic and $(V,pV)$ is a coefficient ring of the complete local ring $R$. So the statement follows from the argument in the proof of the preceding part.
        \end{proof}
    \end{cor}
    
    \section{Reduction  to  excellent unique factorization domains}
    
    In this section we show that some homological conjectures reduce to  excellent unique factorization domains. Here, by  homological conjectures we mean one of the Small Cohen-Macaulay Conjecture,   Big-Cohen-Macaulay Algebra Conjecture as well as the Monomial Conjecture.
    
    \begin{lem} Suppose that $R$ is a Noetherian local domain which is not a field. Then there exists $a\in \mam$ such that $a$ does not have a square root in $F:=\text{Frac}(R)$.
      \begin{proof}
        We first claim that there exists $a\in R$ without a square root in $\Frac(R)$. Assume to the contrary that each element of $R$ has a square root in $F$. Then $F$ is square root closed. Let $\mathfrak{p}$ be a height one prime ideal of the integral closure $\overline{R}$ of $R$ which is a Krull domain (see, \cite[Theorem 4.10.5]{HunekeSwansonIntegral}).  Since $\overline{R}_\map$ is an integrally closed domain, so the discrete valuation ring $(\overline{R}_\map,x\overline{R}_\map)$ is square root closed in $\Frac(R)$.  But then $x$ has a square root in $\overline{R}_\map$. This observation violates  $x\notin \map^2 \overline{R}_\map$. 
        
        Hence there exists an element $a\in R$ without a square root in $\Frac(R)$. If $a\in \mam$ then we are done, so suppose that $a$ is an invertible element of $R$. Then for any $0\neq t\in \mam$ we have, $t$, or, $ta$, does not have a square root in $\Frac(R)$.
      \end{proof}
    \end{lem}
    
	\begin{rem} \label{NewRing}

	   Suppose that $R$ is a domain. Since the trivial extension $R\ltimes \omega_R$  $R$ by its canonical module (even the amalgamated duplication ring of $R$ with its canonical ideal\footnote{See, \cite{BagheriYassemiAConstruction}.}) is not a domain we will use a different multiplication on $R\dsum \omega_R$. To be more precise, since $R$ is a domain so \cite[(3.1)]{AoyamaGotoOnTheEndomorphism} implies that $\omega_R$ is an ideal of $R$. By applying the preceding lemma we can choose $a\in \mam$ such that $a$ does not have a square root in $\Frac(R)$. We may endow $S:=R\dsum \omega_R$ with a ring structure such that, 
	    $$
	    (r,x)(r',x')=(rr'+xx'a,rx'+r'x),
	    $$
	    for each $(r,x),(r',x')\in R\dsum \omega_R$. In particular, $S$ is a subring of $\Ra$ and it is a domain by Remark \ref{FirstRemark}(v). Furthermore, $\eta_R:R\rightarrow S$ which maps an element $r\in R$ to $(r,0)$ is an $R$-algebra homomorphism. Note that $S$ is a  local ring with the maximal ideal $\man:=\mam\dsum \omega_R$ and $\eta$ is a local monomorphism.  As we will see in the following remark,  $S$ endowed with the aforementioned multiplication is  an integral domain which is a complete intersection in low dimension, provided $R$ is a normal ring. The fact that it is a complete intersection in the low dimension is essential for reducing the homological conjectures to the subclass of unique factorization domains. 
    
     \end{rem}

     \begin{lem} \label{LocallyCompleteIntersection}
       Suppose that $R$ is a normal domain and that $\omega_R$ is the canonical ideal of $R$. Let $a\in \Frac(R)$ be an element of $\mam$ without a square root in $\Frac(R)$. Then the following statements hold.  
       
         \begin{enumerate}
         	\item   [(i)] $R(a^{1/2})$ is locally complete intersection at codimension $\le 1$.
         	\item [(ii)] The subring $R\dsum \omega_R$ of $\Ra$ is locally complete intersection at codimension $\le 1$.
         	         \end{enumerate}
         	  \begin{proof}
         	  	  (i) Let $\map\in \ts\big(R(a^{1/2})\big)$ be a height one prime ideal and set $\maq:=\map\bigcap R$. Firstly note that $\hit_R(\maq)=\hit_{R(a^{1/2})}\map,$ as we have  going down, going up and incomparability.  In particular, $\map\in \text{ass}\big(\maq R(a^{1/2})\big)$. Since $R$ is normal, $R_\maq$ is  a regular local ring. Consequently, $R_\maq[x]$ is a regular domain. Hence $R_\maq[x]/(x^2-a)R_\maq[x]\cong R_\maq(a^{1/2})$ is locally complete intersection. An easy verification shows that,

         	  	  \begin{equation}
         	  	    \label{AssEquality}
         	  	    \unil_{\mathfrak{l}\in \text{ass}\big(\maq \Ra\big)}\mathfrak{l}=
         	  	    \{(r,s)\in \Ra:r^2-s^2a\in \maq\}.
         	  	\end{equation}
         	  	Note that $\maq \Ra=\maq\dsum \maq$. Set, $S:=\Ra\backslash \unil_{\mathfrak{l}\in \text{ass}\big(\maq \Ra\big)}\mathfrak{l}$. Then    there exists the natural embedding $\varphi:R_\maq\rightarrow S^{-1}\big(\Ra\big)$ defined by $a/s\rightarrow (a,0)/(s,0)$. On the other hand $a$ does not have a square root in $\Frac(R_\maq)=\Frac(R)$. Thus by Remark \ref{FirstRemark}(iv) $\varphi$ induces the injection $$\Psi:R_\maq(a^{1/2})\rightarrow S^{-1}\big(R(a^{1/2})\big),$$ such that maps any element $(a/b,a'/b')\in R_\maq(a^{1/2})$ to $(ab',a'b)/(bb',0)$.  For each element  $(r,s)/(r',s')\in S^{-1}\big(R(a^{1/2})\big)$ we have, by(\ref{AssEquality}) , $r'^2-s'^2a\notin \maq$ and,
         	  	    $$ \Psi\big((rr'-ss'a)/(r'^2-s'^2a),(r's-rs')/(r'^2-s'^2a)\big)=(r,s)/(r',s').$$
         	  	Hence, $\Psi$ is an isomorphism. So, $S^{-1}\big(\Ra\big)$, is locally complete intersection. Consequently, as $\map\in \as\big(\maq \Ra\big)$, we have $\Ra_\map\cong S^{-1}\big(\Ra\big)_{S^{-1}\map}$ is a complete intersection.
         	  	    
         	  	 (ii) Let $\map'\in \ts (R\dsum \omega_R)$   be a height one prime ideal. Since $\Ra$ is an integral extension of $R\dsum \omega_R,$ so there exists $\map\in \ts\big(\Ra\big)$ lying over $\map'$. Let $\maq$ be the contraction of $\map'$ in $R$. Note that $\hit_R(\maq)=\hit_{\Ra}(\map),$ because $\Ra$ is an integral flat extension of $R$. On the other hand,  because $R$ is integrally closed so $R\rightarrow R\dsum \omega_R$ has the going down property. Hence, $\hit_{\Ra}(\map)=\hit_{R}(\maq)\le \hit_{R\dsum \omega_R}(\map')$. This shows that $\hit_{\Ra}(\map)=1$, as well. 
         	  	 
         	  	 If $\omega_R\nsubseteq \maq$, then it is easily verified that the natural map $(R\dsum \omega_R)_{\map'}\rightarrow \Ra_\map$ is an isomorphism.
         	  	 
         	  	  Now, we deal with the case where $\omega_R\subseteq \maq$. Note that $(R\dsum \omega_R)_{\map'}$ is a localization of $(R\dsum \omega_R)_\maq$. So, it suffices to show that the latter is (locally) complete intersection. Moreover, $(R\dsum \omega_R)_\maq$ is an $R_\maq$-algebra by the map $r/t\mapsto (r,0)/t$ and this underlying $R_\maq$-module structure of $(R\dsum \omega_R)_\maq$ is the same $R_\maq$ module $(R\dsum \omega_R)_\maq\cong R_\maq\dsum R_\maq$. In particular, $(R\dsum \omega_R)_\maq $ is a flat $R_\maq$-algebra. So it is sufficient to show that the closed fiber of $R_\maq\rightarrow (R\dsum \omega_R)_\maq$ is a complete intersection. Since $(\omega_R)_\maq\subseteq \maq R_\maq$ so it is easily seen that there exists a natural ring isomorphism,
         	  	    $$ (R\dsum \omega_R)_\maq /(\maq \dsum \maq \omega_R)_\maq\rightarrow R_\maq/\maq R_\maq \ltimes (\omega_R)_\maq/\maq(\omega_R)_\maq,$$ such that maps an element $\overline{(r,s)/t}$ to $(\overline{r/t},\overline{s/t})$.  Now  from the fact that $(\omega_R)_\maq\cong R_\maq$, we can deduce  that $R_\maq/\maq R_\maq \ltimes (\omega_R)_\maq/\maq(\omega_R)_\maq$ is isomorphic to to the complete intersection ring, $$(R_\maq/\maq R_\maq)[X]/(X^2).$$
         	  \end{proof}
     \end{lem}
     
     In \cite{UlrichGorenstein} the author proves that any Gorenstein ring which is a homomorphic image of a regular ring and it is a complete intersection at codimension $\le 1$ is a specialization of a unique factorization domain. The first part of the proof of the following proposition is more or less repeating the proof of the Ulrich's result, and is stated here for the convenience of the reader and also for the sake of completeness. 
     
     \begin{prop} \label{QuasiGorensteinDeformsToUFD} Suppose that $P$ is a regular  ring and $R:=P/\maa$ is a quasi-Gorenstein ring. Assume, furthermore, that $R$ is locally complete intersection at codimension $\le 1$. Then there exists a unique factorization domain $S$ (which is of finite type over $P$) and a regular sequence $\mathbf{y}$ of $S$ such that $R\cong S/(\mathbf{y})$.
     	\begin{proof}
     		The general idea, here, is similar as given in \cite[Proposition 1]{UlrichGorenstein}. Firstly, we are going to use the idea of generic linkage presented in      \cite{HunekeUlrichDivisor} to be sure that the almost complete intersection linked to $\maa$ is a prime ideal which is a complete intersection at codimension $\le 1$. To be more precise, let $a_1,\ldots, a_n$ be a system of generators for $\maa$ and assume that $\grd(a)=g$. Set, $Q=P[Y_{i,j}:1\le i\le g, 1\le j\le n]$ and $c_i=\sum\limits_{j=1}^n Y_{i,j}a_j$. Then, in the light of, \cite[Proposition 2.9.(b)]{HunekeUlrichDivisor}  and \cite[Proposition 2.6]{HunekeUlrichDivisor}, the linked ideal $\mab=(c_1,\ldots,c_g):_Q\maa Q$, to $\maa Q$, is a prime almost complete intersection which is a  complete intersection at codimension $\le 1$. In particular, in the light of \cite[(5), page 268]{HunekeOnTheSymmetric} one can deduce that $\mab$ is generated by a d-sequence and thence \cite[Theorem 3.1.]{HunekeOnTheSymmetric} implies that $\text{sym}(\mab)\cong \mathcal{R}(\mab)$ wherein $\mathcal{R}(\mab)$ denotes the Rees algebra of $\mab$ and $\text{sym}(\mab)$ stands for the symmetric algebra of $\mab$. Suppose that, $\mab=(c_1,\ldots,c_g,h)$.  Now, in the light of \cite[Theorem 2.2]{HunekeAlmost} in conjunction with \cite[Thoerem 1]{HochsterCriteria}, the extended Rees algebra $Q[t^{-1},\mab t]$ is a unique factorization domain and whence so is 
     		  $$ U:=Q[t^{-1},\mab t]_{ht}=Q[\mab t]_{ht} =\text{sym}(\mab)_{ht}=(Q[Z_0,\ldots,Z_g]/\mad)_{Z_0},$$ wherein, $\mad=<\{\sum\limits_{i=0}^g r_iZ_i|r_0h+\sum\limits_{i=1}^gr_ic_i=0\}>$. Now,  as stated in \cite[Proposition 1]{UlrichGorenstein}, $U$ has an $\mathbb{Z}$-grading structure such that, $$V:=Q[Z_0,Z_1,\ldots,Z_g]/(\mad+(Z_0-1)),$$ is its degree zero subring and since the unique factorization domain $U$ has an invertible element of degree one it is easily seen that $V$ is also a unique factorization domain.
     		  
     		    The reminder of the proof differs with \cite[Proposition 1]{UlrichGorenstein}. Note that if we localize $V$ at its maximal ideal $\man:=(\mam,Y_1,\ldots,Y_{ng},Z_1,\ldots,Z_g,Z_0-1)$ then $R=V_\man/(Y_1,\ldots,Y_{ng},Z_1,\ldots,Z_g)$. Hence in order to deduce the statement it suffices to show that $Z_1,\ldots,Z_g$ forms a regular sequence of $V$, i.e. $Z_0-1,Z_1,\ldots,Z_g$ is a regular sequence on $\text{sym}(\mab)$. We achieve this with the aid of theory of $\mathcal{Z}$-complexes which is introduced and investigated in   \cite{HerzogSimisKoszul }. We use the notation $Z_i$ to denote the $i$-th cycles of the Koszul complex of Q with respect to the sequence $\mathbf{c}:=h,c_1,\ldots,c_g$. Let $S=Q[Z_0,\ldots,Z_g]$, be the polynomial ring over $Q$. Then the $\mathcal{Z}$-complex  associated to the sequence $\mathbf{c}$ is the complex,
     		      
     		        $$ \mathcal{Z}:=0\rightarrow Z_g\bigotimes_Q S(-g)\rightarrow \cdots\rightarrow Z_2\ten_Q S(-2)\rightarrow Z_1\ten_Q S(-1)\rightarrow S\rightarrow 0, $$
     		        
     		        whose differential is defined by the rule,
     		          $$ \partial \big((\sum\limits_{j} q_je_{j_1}\wedge \cdots \wedge e_{j_i})\bigotimes f\big)=\sum\limits_jq_j\big(\sum\limits_{k=1}^i(-1)^ke_{j_1}\wedge\cdots \wedge \widehat{e_{j_k}}\wedge \cdots \wedge e_{j_i}\bigotimes Z_kf\big).$$
     		          It is easily seen that $\Ht_0(\mathcal{Z})=\text{sym}(\mab)=S/\mad$. Moreover, $\mathcal{Z}$ is acyclic by virtue of  \cite[Theorem 12.5]{HerzogSimisKoszul } because any $d$-sequence is a proper sequence(see, \cite[Definition  6.1.]{HerzogSimisKoszul}). Using these two facts and by considering the spectral sequence arising from the double complex, $$\mathcal{Z}\bigotimes_S \Kt_\bullet(Z_0-1,Z_1,\ldots,Z_g;S),$$ we can deduce that $Z_0-1,Z_1,\ldots,Z_g$ is a regular sequence on $\text{sym}(\mab)$ if and only if, $$\mathcal{Z}\bigotimes_S S/(Z_0-1,Z_1,\ldots,Z_g),$$ remains acyclic. Let $e_0$ denotes the basis of $\Kt_1(h,c_1,\ldots,c_g;Q)=Q^{g+1}$ corresponding to the element $h$. Then each element $\zeta_i$ of $Z_i$ can be represented as $\zeta_i=z_i+z'_i\wedge e_0$ such that $z_i\in \dsum_{j_1\neq 0}  Q_{e_{j_1}\wedge\cdots\wedge e_{j_i}}\subsetneq \Kt_i(h,c_1,\ldots,c_g;Q)$ and $z'_i\in \dsum_{j_1\neq 0}  Q_{e_{j_1}\wedge\cdots\wedge e_{j_{i-1}}}\subsetneq \Kt_{i-1}(h,c_1,\ldots,c_g;Q)$. Note that,
     		            \begin{align*}
     		              \mathcal{Z}\bigotimes_S \big(S/(Z_0-1,Z_1,\ldots,Z_g)\big)= & \\ &
     		                0\rightarrow  Z_g\bigotimes_Q Q[Z_0]/(Z_0-1)\rightarrow \cdots \rightarrow   Z_2\bigotimes Q[Z_0]/(Z_0-1) \rightarrow  & \\ &  \rightarrow Z_1\bigotimes_Q Q[Z_0]/(Z_0-1)\rightarrow Q[Z_0]/(Z_0-1)\rightarrow 0,
     		             \end{align*}
     		             whose differentials satisfy, $\partial\big((z_i+z_i'\wedge e_0)\bigotimes \overline{f}\big)=(-1)^{i-1}z_i'\bigotimes \overline{f},$\footnote{More precisely, note that $\mathcal{Z}\bigotimes_S S/(Z_0-1,Z_1,\ldots,Z_g)$ is a subcomplex of the complex, $$\mathcal{K}:=0\rightarrow \Kt_{g+1}\ten_Q S/(Z_0-1,Z_1,\ldots,Z_g)\rightarrow \Kt_g\ten_Q S/(Z_0-1,Z_1,\ldots,Z_g)\rightarrow \cdots\rightarrow  \Kt_0\ten_Q S/(Z_0-1,Z_1,\ldots,Z_g)\rightarrow 0,$$ wherein $\Kt_i$ denotes the $i$-th component of the Koszul complex $\Kt_\bullet(h,c_1,\ldots,c_g;Q)$. The differential of $\mathcal{K}$  is defined by,
     		             	$ \partial (e_{j_1}\wedge \cdots \wedge e_{j_i}\bigotimes f)=\sum\limits_{k=1}^i(-1)^ke_{j_1}\wedge\cdots \wedge \widehat{e_{j_k}}\wedge \cdots \wedge e_{j_i}\bigotimes Z_kf.$ The  complex $\mathcal{K}$ has a DG-algebra structure. Using this DG-algebra structure of $\mathcal{K}$ the identity can be easily verified because in view of the definition of $\partial$ in conjunction with the fact $Z_i$ has  zero image in $S/(Z_0-1,Z_1,\ldots,Z_g)$ for each $1\le i\le g$ we can conclude that $\partial(z_i)=\partial(z'_i)=0$. This DG-algebra technique is used in the proof of \cite[Lemma 4.3]{HassanzadehNaelitonResidual}.},  for each $z_i+z'_i\wedge e_0\in Z_i$. Let, $\alpha$, be a cycle in $\mathcal{Z}_i$. Without loss of generality we can say, $\alpha =(z_i+z'_i\wedge e_0)\bigotimes \overline{1}$. Thus, $z'_i\bigotimes \overline{1}=0$ which implies that, $z'_i=0$. Hence, $\alpha=z_i\bigotimes \overline{1}$. In particular, $z_i\in Z_i(c_1,\ldots,c_g;Q)$. Now the exactness of $K_\bullet(c_1,\ldots,c_g;Q)$ implies that there exists $\zeta_i\in \Kt_{i+1}(c_1,\ldots,c_g;Q)$ such that $\partial(\zeta_i)=z_i$. In particular, $-h\zeta_i+(-1)^iz_i\wedge e_0\in Z_{i+1}$. Now using the algebra structure we get,
     		               $$ \alpha=z_i\bigotimes \overline{1}=\partial\bigg(\big(-h\zeta_i+(-1)^iz_i\wedge e_0\big)\bigotimes \overline{1}\bigg). $$ This observation concludes the proof.
     		                   
     	\end{proof}     	
     \end{prop}
      Recall that by  homological conjectures we mean one of the Small Cohen-Macaulay Conjecture,   Big-Cohen-Macaulay Algebra Conjecture as well as the Monomial Conjecture. The first part of the subsequent theorem reduces the homological conjectures 	to the class of quasi-Gorenstein domains. As pointed out in the introduction, as the main result of this section, the second part of the subsequent theorem gives the promised descent to excellent unique factorization domains. The integral closure of a complete quasi-Gorenstein ring is not necessarily quasi-Gorenstein. Therefore,  we could not reduce to the normal complete quasi-Gorenstein domains directly by taking integral closure in the first part. But by virtue of the second part this is achievable which is stated separately in the third part. Note that, as stated in \cite{UlrichGorenstein}, the completion of the unique factorization domains in  \cite{UlrichGorenstein}  are not unique factorization domain necessarily. However any possible reduction of the homological conjectures to the complete unique factorization domains would be fairly helpful as the complete unique factorization domains have 	even more interesting properties at least in the equal characteristic zero. More precisely, any equal characteristic zero complete unique factorization domain with algebraically closed residue field of dimension up to $4$ is Cohen-Macaulay.  Furthermore, with mild conditions, any complete unique factorization domain of equal characteristic zero with algebraically closed residue field  satisfies the Serre-condition $S_3$ (see, \cite[page 540]{LipmanUnique} or \cite{HartshorneOgusOnTheFactoriality} for both of the aforementioned results).

     \begin{thm} \label{mainResult} The following statements holds.
     	\begin{enumerate}
     	\item[(i)]
     	If each complete quasi-Gorenstein  domain, which is locally complete intersection at codimension $\le 1$, satisfies the  homological  conjectures then, each  local ring  satisfies homological conjectures. This reduction is dimensionwise, in the sense that  homological conjectures for $d$-dimensional local rings reduce to the $d$-dimensional complete quasi-Gorenstein domains.
     	\item [(ii)]
     	The homological conjectures reduce to the excellent (and homomorphic image of regular) unique factorization domains. This reduction is dimensionwise for the Monomial Conjecture. 
     	\item[(iii)]  
     		The homological conjectures reduce to quasi-Gorenstein normal complete domains. This reduction is dimensionwise for the Monomial Conjecture.
       \item[(iv)] For the Canonical Element Conjecture it suffices only to consider mixed characteristic complete normal quasi-Gorenstein rings of minimum possible depth, i.e. depth $2$.
     \end{enumerate}
     \begin{proof}
     	  (i) We may suppose that $R$ is a complete normal domain. In particular $R$ satisfies the $S_2$-condition. Let $\omega_R$ be the canonical ideal of $R$ and set $S:=R\dsum \omega_R$ as stated in the Remark \ref{NewRing}. Note that $S$ is a complete local ring. By \cite[Proposition 1.2]{AoyamaGotoOnTheEndomorphism} we have,
     	  
     	  \begin{align*}
     	  	\homm_{R}\big(H_{\mathfrak{n}}^{d}(S),E(R/\mathfrak{m})\big)&\cong \homm_{R}\big(H_{\mathfrak{m}}^{d}(\omega_{R})\bigoplus H_{\mathfrak{m}}^{d}(R),E(R/\mathfrak{m})\big) &\\ & \cong \homm_{R}\big(H_{\mathfrak{m}}^{d}(R)\bigotimes_{R}\omega_{R},E(R/\mathfrak{m})\big)\bigoplus \homm_{R}\big(H_{\mathfrak{m}}^{d}(R),E(R/\mathfrak{m})\big) &\\& \cong \homm_{R}(\omega_{R},\omega_{R})\bigoplus\omega_{R}& \\ & \cong R\bigoplus\omega_{R}=S.
     	  \end{align*}
     	  Hence  we get an isomorphism, $S\rightarrow \homm_R\big(\Ht^d_\man(S),\Et(R/\mam)\big)$, of $R$-modules which is an $S$-isomorphism too\footnote{In order to see why this is  also an $S$-isomorphism, we fix an $R$-isomorphism, $$\varphi:\omega_R\rightarrow \homm_R\big(\Ht^d_\mam(R),\Et(R/\mam)\big).$$ Then the composition of the above isomorphisms which is the mentioned $R$-isomorphism $S\rightarrow \homm_R\big(\Ht^d_\man(S),\Et(R/\mam)\big)$ maps an element $(r,\alpha)\in R\dsum \omega_R=S$ to the $R$-homomorphism defined by the rule, $$[(r',\alpha')+(\mathbf{x}^n)]\in \Ht^d_\man(R\dsum \omega_R)\mapsto \big(\varphi(r\alpha'+r'\alpha)\big)\big([1+(\mathbf{x}^n)]\big),$$ wherein $\mathbf{x}$ is a system of parameters for $R$. Then it is easily verified that the map is, moreover, an $S$-homomorphism.}. In particular, $H^d_{\mathfrak{n}}(S)\cong \homm_{R}(S,E(R/\mathfrak{m}))$, as $S$-modules because on the category of $S$-modules two functors $\homm_R(-,\Et(R/\mam))$ and  $\homm_S(-,\Et(S/\man))$ are naturally equivalent . So,
     	  \begin{align*}
     	  	\homm_{S}\big(S/\mathfrak{n},H_{\mathfrak{n}}^{d}(S)\big)\cong \homm_{S}\Big(S/\mathfrak{n},\homm_{R}\big(S,E(R/\mathfrak{m})\big)\Big)& \cong \homm_{R}\big(S/\mathfrak{n}\bigotimes_{S}S,E(R/\mathfrak{m})\big)& \\ & \cong \homm_{R}\big(R/\mathfrak{m},E(R/\mathfrak{m})\big)& \\ &\cong R/\mathfrak{m}& \\ &\cong S/\mathfrak{n}. 
     	  \end{align*}
     	  This observation implies that $\Ht^d_\man(S)$ has one dimensional socle. Since $a$ does not have a square root in $\Frac(R),$  so Remark \ref{FirstRemark} (v) implies that $R(a^{1/2})$ and whence its subring $S$ is an integral domain. Thus, by \cite[Theorem 5.7]{TavanfarTousiAStudy}, $S$ is a quasi-Gorenstein domain. Since $S$ is a finitely generated (as $R$-module) complete quasi-Gorenstein extension domain of $R$ so the statement follows.
     	  
     	  (ii) Using  the previous part in conjunction with Lemma \ref{LocallyCompleteIntersection} and Proposition  \ref{QuasiGorensteinDeformsToUFD} we may assume that $R$ is a complete quasi-Gorenstein domain such that $R=S/\mathbf{y}S$ wherein $(S,\man)$ is a unique factorization domain and $\mathbf{y}:=y_1,\ldots,y_n$ is a regular sequence contained in $S$. Clearly, if $S$ satisfies the homological conjectures then  so does $R$. This proves the first claim. For proving the second claim  we invoke the Hochster's  first general grade reduction technique developed in \cite{HochsterProperties}.  Let $a_1,\ldots,a_d$ be a system of parameters for $R$. Assume that the Monomial Conjecture is valid for every $d$-dimensional unique factorization domain but, to the contrary, $$a_1^t\cdots a_d^t\in (a_1^{t+1},\ldots,a_d^{t+1}),$$ for some $t\in \mathbb{N}$. So we get the system of parameters $y_1,\ldots,y_n,y_{n+1},\ldots,y_{n+d}$ for $S$ such that  $y_{i+n}$ is a lift of $a_i$ in $S$ for each $1\le i\le d$ and $$y_{n+1}^t\cdots y_{n+d}^t\in (y_1,\ldots,y_n,y_{n+1}^{t+1},\ldots,y_{n+d}^{t+1}).$$ Since $R$ satisfies the $S_2$-condition so $y_1,\ldots,y_n,y_{n+1}^{t+1},y_{n+2}^{t+1}$ is a regular sequence of $S$.  As $n+2\ge 3$ so \cite[Theorem. c)]{HochsterProperties} implies that the first general grade reduction, $$T=S[X_1,\ldots,X_{n+2}]/(\sum\limits_{k=1}^{n}y_kX_k+y_{n+1}^{t+1}X_{n+1}+y_{n+2}^{t+1}X_{n+2}),$$ of  $S$ is again a (non-local) unique factorization domain and so is its localization at $\man T$. We have, $\dime(T_{\man T})=\dime(S)-1$. Furthermore, as $X_1$ is invertible in $T_{\man T}$ so $y_2,\ldots,y_{n+d}$ is a system of parameters for $T_{\man T}$, $y_2,\ldots,y_{n},y_{n+1}^{t+1},y_{n+2}^{t+1}$ is a regular sequence of $T_{\man T}$ and $$y_{n+1}^t\cdots y_{n+d}^t\in (y_2,\ldots,y_n,y_{n+1}^{t+1},\ldots,y_{n+d}^{t+1})T_{\man T}.$$ Proceeding in this way, we get a unique factorization domain $U$ with $\dime(U)=\dime(R)=d$ such that the image of $y_{n+1},\ldots,y_{n+d}$ in $U$ is a system of parameters for $U$ which does not satisfy the Monomial Conjecture. This is a contradiction.
     	  
     	  (iii) This is followed by taking the completion of the unique factorizations of the second part, because by \cite[Lemma (2.4)]{FossumFoxbyMinimal} any unique factorization domain with a canonical module is quasi-Gorenstein. Note also that the completion remains normal, as the formal fibers of an excellent ring are regular.
     	  
     	  (iv)   In order to see why this is the case firstly note that by \cite{HochsterCanonical} the Canonical Element Conjecture for all local rings is valid if and only if all local rings satisfies the Monomial Conjecture. Therefore in view of   Theorem \ref{mainResult} (iii) the Canonical Element Conjecture reduces to the mixed characteristic normal quasi-Gorenstein complete domains. So let $R$ be such a quasi-Gorenstein normal domain. Without loss of generality we may assume that $R$ has infinite residue field. If $\dep(R)\gneq 2$ then applying \cite[Theorem 4.4(Local Bertini Theorem)]{OchiaiShimomotoBertini} there exist $x\in R$ such that $R/xR$ is a normal domain. In particular $R/xR$ satisfies the $S_2$-condition and thence $R/xR$ is a quasi-Gorenstein normal domain in view of \cite[Corollary 3.4.(i)]{TavanfarTousiAStudy}. Consequently, \cite[Proposition 3.1.]{TavanfarTousiAStudy} implies that $x\notin \unil_{\map\in \text{Att}\big(\Ht^{d-1}_\mam(R)\big)}\map$. In particular, if $R/xR$ satisfies the Canonical Element Conjecture then so does $R$ in the light of \cite[(5.3) Theorem]{HochsterCanonical}. Proceeding in this fashion the claim follows.   
     \end{proof}
     \end{thm}
     
     The theory of generic linkage had an important role in the reduction  of the homological conjectures to unique factorization domains. In the second part of the subsequent remark we, again, observe that how the linkage theory is related to the homological conjectures.  In the first part we show that the validity of the Canonical Element Conjecture in an open subset implies the validity in general.
     
     \begin{rem} \label{ReductionRemark}
     	(i) The following approach to the Monomial Conjecture might be helpful: The Canonical Element Conjecture is valid if it is valid, in some sense, in an open dense subset. More precisely, in view of the first part of  Theorem \ref{mainResult}(iii), we can assume that $R$ is a complete normal quasi-Gorenstein domain of mixed characteristic with $\dep(R)\ge 3$.  Let $(V,\pi_V,k)$ be a  coefficient ring of $R$.  Let $x_0,\ldots,x_n$ be a minimal generating set of $\mam$. If there exists  a dense Zariski open subset  $O\subseteq \mathbb{P}^d(k)$ such that for any $a=(a_0:\cdots:a_n)$  in the preimage of $O$ in $\mathbb{P}^n(V)$, as described in \cite[Definition 2.1.]{OchiaiShimomotoBertini}, we have $R/\sum\limits_{i=0}^na_ix_i$ satisfies the Canonical Element Conjecture then the Canonical Element Conjecture holds. Indeed, if this is the case then we can apply \cite[Theorem 4.4(Local Bertini Theorem)]{OchiaiShimomotoBertini}  to find an element $a\in R$ such that $R/aR$ is a normal domain satisfying the Canonical Element Conjecture.  Then \cite[Theorem 3.4(i)]{TavanfarTousiAStudy} and \cite[Propositino 3.1.]{TavanfarTousiAStudy} in conjunction with \cite[(5.3) Theorem]{HochsterCanonical} imply that $R$ satisfies the Canonical Element Conjecture.
     	
     	(ii) In \cite[page 158]{StrookerStuckradMonomial} the following question is proposed. Suppose $\maa$ and $\mab$ are ideals in a $d$-dimensional Gorenstein ring $R$ which are geometrically linked w.r.t. the empty regular sequence. Suppose $R$ possess a maximal Cohen-Macaulay module $M$ with $0:_RM=\maa$. Does $R$ possess a maximal Cohen-Macaulay module $N$ with $0:_RN=\mab$? We show that if this question has an affirmative answer then Small Cohen-Macaulay Conjecture is true in dimension $3$. Furthermore if Monomial Conjecture analogous of this question has an affirmative answer then the Monomial Conjecture holds in any dimension. 
     	
     	By Theorem\ref{mainResult}(i) the Small Cohen-Macaulay Conjecture reduces, dimensionwise, to the class of complete quasi-Gorenstein domains. Let $P$ be a regular local ring and $\maa$ be a quasi-Gorenstein prime ideal of $P$. As in the proof of \cite[1.2. Proposition]{DuttaTheMonomial} Choose a maximal regular sequence $x_1,\ldots,x_h$ contained in $\maa$ such that $\maa P_\maa=(x_1,\ldots,x_h)P_\maa$.  Then, $\maa$ and $(x_1,\ldots,x_h):_P\maa$ are geometrically linked and since $\maa$ is a quasi-Gorenstein ideal so $\mab:=(x_1,\ldots,x_h):_P\maa$ is an unmixed almost complete intersection. Therefore $P/\mab$ has positive depth and whence Lemma \ref{CanonicalModule}(ii) implies that $\dep(\omega_{P/\mab})\ge 3$. In particular, in the case of dimension three we have $\omega_{P/\mab}$ is a maximal Cohen-Macaulay $P/\mab$-module  and therefore an affirmative answer to the mentioned question implies the validity of the Small Cohen-Macaulay Conjecture in dimension $3$. 
     	
     	For the Monomial Conjecture analogous of the Strooker and St\"uckrad's question assume in addition that $\maa$ is a normal quasi-Gorenstein prime ideal of $P$. Then
     	the claim follows by \cite[step 2., page 238]{DuttaGriffithIntersection}. In spite of this observation,  it is by no means clear for us that how the theory of linkage can be helpful for this type of problems, in this point of view.      

     \end{rem}

    Any successful proof of the Monomial Conjecture for the class of generalized Cohen-Macaulay domains will settle the Monomial Conjecture in dimension three (by means of a new proof). Although we do not know how to prove the Monomial Conjecture for this type of rings but we show that again the problem reduces to the class of quasi-Gorenstein generalized Cohen-Macaulay rings.

    \begin{prop} \label{quasi-Gorenstein_MC} The following statements hold.
    	\begin{enumerate}
    		
    		\item[(i)]  If every complete quasi-Gorenstein generalized Cohen-Macaulay ring satisfies the Monomial Conjecture then so does every generalized Cohen-Macaulay ring.
    		\item[(ii)] Monomial Conjecture holds for every complete  generalized Cohen-Macaulay domain if and only if every  complete generalized Cohen-Macaualy quasi-Gorenstein domain satisfies Monomial Conjecture.
    	\end{enumerate}  
    	\begin{proof}

    		(i) Let $R$ be a complete generalized Cohen-Macaulay ring. Denote by $S$ the $S_2$-ification of $R$. By \cite[Corollary 4.3]{AoyamaSomeBasic} $\omega_{R_\map}$, for each $\map\in \ts(R)$, is either the canonical module of the Cohen-Macaulay ring $R_\map$ or is zero. Hence, for every $\map\in \ts(R)\backslash \mam$ either $S_\map$ is zero or $S_\map\cong \homm_{R_\map}\big((\omega_R)_\map,(\omega_R)_\map\big)\cong \homm_{R_\map}(\omega_{R_\map},\omega_{R_\map})\cong R_\map$. On the other hand we have, $$\As_R(S)=\As_R\big(\homm_R(\omega_R,\omega_R)\big)=\text{Supp}_R(\omega_R)\bigcap \As_R(\omega_R)=\text{Assh}_R(R).$$ Therefore,  by \cite[9.5.7 Exercise]{BrodmannSharpLocalCohomology} $S$ is a generalized Cohen-Macaulay $R$-module. Let $\man\in \text{Max}(S)$ such that $\dime(S)=\dime(S_\man)$. Now, we show that $S_\man$ is a generalized Cohen-Macaulay ring. Denote the natural map $R\rightarrow S$ by $\eta$. We know that $S$ is a finitely generated $R$-module. So $\dime(R/\text{ker}\ \eta)=\dime(S)$. By \cite[(4.2) Lemma]{HochsterCanonical} we have $a:=\text{ker}(\eta)\subseteq \unil_{\map\in \text{Assht}(R)}\map$ and whence $\dime(S)=\dime(R)$.  Since $S$ is integral over $R/\maa$ so $\mam/\maa=\man\bigcap R/\maa$ and thereby $\mam S_{\man}$ is an $\man S_\man $-primary ideal by the incomparability property of $R/\maa \rightarrow S$. There exists $l\ge 0$ such that $\mam^l\Ht^i_\mam(S)=0$ ($0\le i\le \dime(S)-1$). So that $\mam^l\Ht^i_{\mam S_\man}(S_\man)=0$ ($0\le i\le \dime(S)-1$) for some $l\ge 0$. Thus, $(\man S_\man)^k\Ht^i_{\man S_\man}(S_\man)=0$ ($0\le i\le \dime(S_n)-1$). Hence $S_\man$ is a generalized Cohen-Macaulay ring. We have $S_\man$ is $S_2$. Hence by passing to the completion if necessary, we can assume that $R$ is an $S_2$-complete generalized Cohen-Macaulay ring. Therefore, in view of \cite[Theorem 2.11.]{AoyamaSomeBasic} the trivial extension $T$ of $R$ by its canonical module is a complete quasi-Gorenstein ring. We prove that $T$ is a generalized Cohen-Macaulay ring. Clearly if $T$ satisfies the Monomial Conjecture then so does $R$.
    		
    		Let $\maq\in \ts(T)\backslash \{\mam\ltimes \omega\} $. Then there exists $\map\in \ts(R)\backslash \{\mam\}$ such that $\maq=\map\ltimes \omega_R$. So $T_\maq\cong (R\ltimes\omega_R)_{\map\ltimes \omega_R}\cong R_\map\ltimes (\omega_R)_\map \cong R_\map\ltimes \omega_{ R_\map }$. Hence, $T_\maq$ is a Cohen-Macaulay ring. As $T$ is equidimensional, \cite[9.5.7 Exercise]{BrodmannSharpLocalCohomology} implies that $T$ is generalized Cohen-Macaulay.
    		
    		(ii) Let $R$ be a complete generalized Cohen-Macaulay domain. Set, $S:=\homm_R(\omega_R,\omega_R)$. Then, by virtue of \cite[(3.6)]{HochsterHunekeIndecomposable} $S$ is local\footnote{Since $S$ is local so we do not need to localize $S$. This is important because the completion of  a domain is not necessarily a domain.}. Also we know that $S$ is complete and domain. By a similar argument as in the preceding part, $S$ is generalized Cohen-Macaulay. Hence, as stated in the proof of Theorem \ref{mainResult}(i), $T=S\dsum \omega_S\subseteq S(a^{1/2})$, is a local quasi-Gorenstein complete domain wherein  $a$ is an element of $S$ without a square root in $\Frac(S)$. It is enough to show that  $T$ is a generalized Cohen-Macaulay ring. 
    		
    		Let $\map\in \ts(S)\backslash\{\mam\}$. Then $(S\dsum \omega_S)_\map\cong S_\map\dsum (\omega_S)_\map\cong S_\map\dsum \omega_{S_\map},$ is a Cohen-Macaulay $S_\map$-module. On the other hand, $\As_S(S\dsum \omega_S)=\As_S(S)\uni \text{Assh}_S(\omega_S)=\As_S(S)=\{0\}$. Hence $T$ is a generalized Cohen-Macaulay $S$-module. Now, a similar argument as in the preceding part shows that  $T$ is a generalized Cohen-Macaulay domain. 
    	\end{proof}    	
    \end{prop}

    \section{A class of rings satisfying the Monomial Conjecture}
    
    As stated  in the introduction for the validity of the Monomial Conjecture it suffices only to check the almost complete intersections and system of parameters of the form $p,x_2,\ldots,x_d$ where $p$ is the residual characteristic. This section is devoted to prove either the Small Cohen-Macaulay Conjecture or the Big Cohen-Macaulay Algebra Conjecture for an almost complete intersection $R$ satisfying the extra assumption that,
      \begin{equation}
        \label{StrongHypothesis}
        \mathfrak{m}^2\subseteq (\mathbf{x})R,
      \end{equation}  
      for some system of parameters $x_1,\ldots,x_d$ of $R$. In the mixed characteristic case we moreover assume that $x_1=p$.
      
       A motivation behind this consideration is that this case seems to be the most simple non-trivial case of the Monomial Conjecture. For another motivation see the introduction. But, parallel, to this aim, several other interesting facts about this class of almost complete intersections are also given. The following lemma is not required for the next theorem which is the main result of the section but it is needed for its subsequent remark.
      \begin{lem} \label{CanonicalModule}
      	Assume that $(A,\man)$ is a regular local ring. Let $\maa=(x_1,\ldots,x_d)$ (resp.  $\mab=(y_1,\ldots,y_s)$) be a complete intersection  (resp. an almost complete intersection) ideal of $A$.  Suppose furthermore that,  $x_1,\ldots,x_d,y_1,\ldots,y_{s-1}$, forms a regular sequence of $A$ and that $\hit(\maa)+\hit(\mab)=d+s-1=\dime(A)$. Set, $R:=A/\mab$. Then the following statements hold.
      	\begin{enumerate}
      		\item[(i)] There exists an exact sequence,
      		\begin{center}
      			$0\rightarrow \Ht_2(\mathbf{x},R)\rightarrow \omega_R/\maa \omega_R\rightarrow \omega_{R/\maa}\rightarrow \Ht_1(\mathbf{x},R)\rightarrow 0$.
      		\end{center}
      		
      		\item[(ii)] We have $\Ht_i(\mathbf{x},R)\cong \Ht_{i-2}(\mathbf{x},\omega_R)$ for each $i\ge 3$. In particular, $$\begin{cases}\dep(\omega_R)=\dep(R)+2,\ & \dep(R)\le d-2\\ \dep(\omega_R)=d,\ & \dep(R)\ge d-1.\end{cases}$$
      		
      		\item[(iii)] The system of parameters $\mathbf{x}$ of $R$ satisfies the Monomial Conjecture if and only if, 
      		  $$ \ell\big(\Ht_2(\mathbf{x},R)\big)\lneq \ell\big(\omega_R/\maa \omega_R\big) .$$ Definitely, this is equivalent to say that, $\omega_R\rightarrow \omega_{R/\maa},$ is non-zero in the above exact sequence.

      	\end{enumerate}     
      	\begin{proof}
      	  (i) and (ii):	We have, $$\Ht_1(x_1,\ldots,x_d;R)\cong \Ht_i\Big(K_\bullet\big(x_1,\ldots,x_d;A/(y_1,\ldots,y_{s-1})\big)\ten_A (A/y_sA)\Big).$$ 
      		This encouraged us to consider the spectral sequence, $\Kt_\bullet\big(\mathbf{x};A/(\mathbf{y}')\big)\ten_A \Kt_\bullet\big(y_s;A\big)$, where $\mathbf{y}'$ denotes the truncated sequence $y_1,\ldots,y_{s-1}$. So we take into account the bicomplex $M_{p,q}:=\Kt_p\big(\mathbf{x};A/(\mathbf{y}')\big)\ten_A \Kt_q\big(y_s;A\big)$ in which here, as usual, $p$ stands for the column $p$. Note that, \begin{equation}
      		  \label{SpectralSequnece}
      		\begin{cases}
      		  \Ht_i\big(\text{Tot}(M)\big)\cong \Ht_i\big(\mathbf{x},y_s;A/(\mathbf{y}')\big)\cong \Ht_i\big(y_s;A/(\mathbf{x},\mathbf{y}')\big)=0, &  i\ge 2. \\	       		    \Ht_i\big(\text{Tot}(M)\big)\cong \Ht_i\big(\mathbf{x},y_s;A/(\mathbf{y}')\big)\cong \Ht_i\big(y_s;A/(\mathbf{x},\mathbf{y}')\big)\cong \omega_{R/\maa}, & i=1. 
      		\end{cases}
      		\end{equation}

      		 Furthermore we have, $$^{II}\Et^2_{p,q}=\begin{cases}\Ht_q(\mathbf{x};R),&p=0\\ \Ht_q\Bigg(\mathbf{x};\bigg(\big((\mathbf{y}'):y_s\big)/(\mathbf{y}')\bigg)\Bigg)=\Ht_q(\mathbf{x};\omega_R),& p=1\\0, &p\neq 0,1.\end{cases}$$ Now the desired exact sequence is just the five term exact sequence of this spectral sequence (see, \cite[Theorem 10.31 (Homology of Five-Term Exact Sequence)]{Rotman}).
      		
      		For the second part note that according to the vanishings of (\ref{SpectralSequnece}) for $i\ge 2$, all of the maps, $$ d^2:\Ht_{i+2}(\mathbf{x};R)\rightarrow \Ht_i(\mathbf{x};\omega_R),\ (i\ge 1)$$ arising from the second page of the spectral sequence are isomorphisms.
      		  
      		  (iii): By \cite[1.3. Proposition]{DuttaTheMonomial} our statement is equivalent to the assertion that $\ell\big(\Ht_1(\mathbf{x};R)\big)\lneq \ell(R/\maa)$. Hence in the light of the exact sequence of the first part we are done once we can show that, $$\ell(R/\maa)=\ell(\omega_{R/\maa})=\ell\bigg(\big((\mathbf{x},\mathbf{y}'):y_s\big)/(\mathbf{x},\mathbf{y}')\bigg).$$
      		     But this is clear from the following exact sequence,
      		       $$0\rightarrow \big((\mathbf{x},\mathbf{y}'):y_s\big)/(\mathbf{x},\mathbf{y}') \rightarrow R/(\mathbf{x},\mathbf{y}')\overset{y_s}{\rightarrow} R/(\mathbf{x},\mathbf{y}')\rightarrow R/\maa\rightarrow 0.$$
      	\end{proof}
      \end{lem}
      
      In the  following example we consider a subclass of almost complete intersections satisfying (\ref{StrongHypothesis}) and we show that they satisfy the Small Cohen-Macaulay Conjecture. The main idea behind this observation is that this class of almost complete intersections are the same free square root extensions $\Ra$ of Remark \ref{FirstRemark}. This is not the case for the more general setting of Theorem \ref{MonomialConjectureHolds}. However, in Theorem \ref{MonomialConjectureHolds},  with aid of the multiplicity theory we will again descend to a similar case of quadratic extensions.
       
      \begin{exam}
      	 Suppose that $R=(V,pV)[[X_2,\ldots,X_d,Z_1,Z_2]]/I$, such that $I=(g_1-Z_1^2,g_2-Z_2^2,g_3-Z_1Z_2)$ wherein $g_1,g_2,g_3\in V[[X_2,\ldots,X_d]]$ are non-zero. Set, $(A,\man):=V[[X_2,\ldots,X_d,Z_1,Z_2]]$. Presume that $R$ is an almost complete intersection, i.e. $\hit(I)=2$.   Note that,
      	   $$ p,x_2:=X_2+I,\ldots,x_d:=X_d+I,$$
      	 forms a system of parameters for $R$ satisfying (\ref{StrongHypothesis}). In this example we show that $R$ has a maximal Cohen-Macaulay module. Setting,
      	   $$ S:=V[[X_2,\ldots,X_d,Z_1]]/(g_1-Z_1^2),$$  then $S$ is a free extension of $B:=V[[X_2,\ldots,X_d]]$ wherein $Z_1$ is the square root of $g_1$, i.e. $S=B(g_1^{1/2})$. If $g_1$ has a square root in $B$ then the polynomials $g_1-Z_1^2$ is reducible and therefore any $\map\in \text{assht}(I)$ contains an element $a\in \man\backslash \man^2$. In particular, then, $\map/aA$ is a height one prime ideal of the regular local ring $A/aA$ which implies that $\map$ is generated by two elements and thence $A/\map$ is a complete intersection and a maximal Cohen-Macaulay $R$-module. So without loss of generality we can assume that $g_1$ does not have a square root in $B$ and whence in $\Frac(B)$.  Consequently, $S$ is a domain  hypersurface by Remark \ref{FirstRemark}(v). Next set, $$ T:=V[[X_2,\ldots,X_d,Z_1,Z_2]] /(g_1-Z_1^2,g_2-Z_2^2).$$   Therefore, $T$ is the free extension of $S$ wherein $Z_2$ is a square root of $g_2$, i.e. $T=S(g_2^{1/2})$. In particular, $T$ is a domain if and only if $g_2$ does not have a square root $\Frac(S)$. But since our almost complete intersection $I$ is an ideal of height two generated by three elements so indeed $g_2$ has a square root in $\Frac(S)$. Therefore by Remark \ref{FirstRemark}(iv) there exists an induced map from $T$ to the integral closure $\overline{S}$ of $S$, denoted by $\Psi$. Since the kernel of this map is an associated prime of the defining ideal of $T$ so, as will be seen later, without loss of generality we can presume that the kernel of $\Psi$ contains the ideal, $I$, and thence we get an induced map $R\rightarrow \overline{S}$. Hence it is enough to prove the statement for $\overline{S}$ (recall that $\overline{S}$ is finitely generated as $S$-module). Note that, $$[\Frac(\overline{S}):\Frac(B)]=[\Frac(S):\Frac(B)]=2.$$ Here the residual characteristic of $B$ is either two or $p>2$. In the first case \cite[Theorem 3.8.]{KatzOnTheExistence} implies that $R$ has a maximal Cohen-Macaulay module, while in the second case $\overline{S}$ is a Cohen-Macaulay ring in the light of \cite{RobertsAbelian}\footnote{By an argument as in  the  proof of Theorem \ref{MonomialConjectureHolds} we can deduce that $\overline{S}$ is necessarily Cohen-Macaulay regardless of the residual characteristic, because we have a quadratic extension. This fact is asserted in Theorem \ref{CohenMacaulayQuadratic}.}. 
      	     
      	      It remains to show that why we can assume that, $I\subseteq \text{ker}\ \Psi$. We claim that $T$  has at most two associated primes.  In order to see why this is the case note that, in the light of \cite[Theorem 1.]{FoxbyInjective}, $E\ten_S T$ is an injective $T$-module provided $E$ is an injective $S$-module. In particular, $\text{E}(S)\ten_S T$, is an injective $T$-module which contains $T$ and thence contains $\text{E}_T(T)$.  Therefore, $\text{E}_T(T)$, as $S$-module, is an injective submodule of $\text{E}(S)\dsum \text{E}(S)$. Since $S$ is a domain so $\text{E}(S)$ is an indecomposable injective module. This observation shows that $\text{E}_{T}(T)$ has at most two direct summands and, as an immediate consequence, $T$ has at most two associated primes. Now if $\alpha_1$ and $\alpha_2$ are distinct square roots of $g_2$ in $\Frac(S)$\footnote{Recall that $S$ has mixed characteristic and therefore $\Frac(S)$ has characteristic zero. So we have two distinct roots.} then it is easily verified that  the rules $(r,s)\mapsto r+s\alpha_i,\  (1\le i\le 2)$ define distinct ring homomorphisms, $\Psi_i:T\rightarrow \overline{S}$, with distinct kernels, $\map_1,\map_2\in \As(T)$. Now the claim, easily, follows from the fact that $T$ has at most two associated primes.
      \end{exam}
      
      \begin{exam}
      	It is worth to give a concrete and explicit non-Cohen-Macaulay example of the class of 	the rings of the previous example. We do this, but in the equal characteristic case, where we can use the Macaulay2 system for computations. Set,  
      	    $$ R=K[[Y_1,\ldots,Y_6,Z_1,Z_2]]/I,$$
      	    such that,
      	      $$ I=(Y_2^6Y_3^5+Z_2^2,Y_3^3Y_4^8+Z_1^2,Y_2^3Y_3^4Y_4^4+Z_1Z_2).$$
      	   Then we have, $\dep(R)=4$  and $\dime(R)=6$. In particular, Lemma \ref{CanonicalModule} implies that $\omega_R$ is, also, a maximal Cohen-Macaulay $R$-module.
      \end{exam}
      
      \begin{thm}\label{MonomialConjectureHolds}
      	Suppose that $R$ is an almost complete intersection satisfying (\ref{StrongHypothesis}). Then $R$ has a big Cohen-Macaulay algebra. A fortiori, $\widehat{R}$ has a maximal Cohen-Macaulay module which is an $\widehat{R}$-algebra.
      	  \begin{proof}
      	  	If $R$ is a Cohen-Macaulay then there is nothing to prove, so assume that $R$ is not Cohen-Macaulay (and is complete). The multiplicity theory plays a pivotal role in the proof. The same idea of the proof of Corollary \ref{WeCanAssumeMinimalGenerator} can be applied to impose the extra assumption that $\mathbf{x}$ is a part of minimal generating set of $\mathfrak{m}$\footnote{Here it is important to point out that the extensions used in the proof of Corollary \ref{WeCanAssumeMinimalGenerator} preserves the condition (\ref{StrongHypothesis}) as well as the almost complete intersection property.}.  Consider a minimal Cohen-presentation $R=A/I$ of $R$ so  $\text{embdim}\big((A,\man)\big)=\text{embdim}\big((R,\mathfrak{m})\big)$ and $I\subseteq \man^2$. Firstly, we have,
      	  	  $$ \ell\big(R/(\mathbf{x})\big) =\ell\bigg(A/\big((\mathbf{x})+I)\big)\bigg)=\ell(\overline{A}/\overline{I}),$$ wherein the notation, $\overline{\ \ }$, means modulo $(\mathbf{x})$. But, $\overline{I}=\overline{\man^2},$ as $\man^2\subseteq I+(\mathbf{x})$ and $I\subseteq \man^2$ simultaneously. It turns out that, 
      	  	  
      	  	     $$ \ell\big(R/(\mathbf{x})\big) = \ell (\overline{A}/\overline{\man}^2)=\text{embdim}(\overline{A})+1=\text{embdim}(R)-\dime(R)+1.$$
      	  	   We claim that, $\text{embdim}(R)-\dime(R)\le 2$, but before proving this claim let us conclude the statement from it. The above inequality in conjunction with our claim implies that $\ell\big(R/(\mathbf{x})\big)\le 3$.    Now, in the light of \cite[Theorem 1, page. 57]{SerreLocalAlgebra}  together with \cite[Corollary , page 90]{SerreLocalAlgebra} we can deduce  a very strong implication that $\et(\mathbf{x},R)\le 3$. But if $\et(\mathbf{x},R)= 3$ then  $R$ is Cohen-Macaulay by \cite[(1.5) Proposition]{HermannIkedaEquimultiplicity}. Thus we have even stronger condition that $\et(\mathbf{x},R)\le 2$. In the case where $\et(\mathbf{x},R)=1$ using \cite[(1.10.1)]{HermannIkedaEquimultiplicity} and \cite[(6.8) Theorem]{HermannIkedaEquimultiplicity} there exists a prime ideal $\map\in \text{Assht}(R)$ such that $R/\map$ is regular, hence is s a  maximal Cohen-Macaulay $R$-module. It remains to deal with the case where $\et(\mathbf{x},R)=2$. In this case by a similar argument as above without loss of generality we can assume that there exists $\map\in \text{Assht}(R)$ for which $R/\map$ satisfies $\et(\mathbf{x},R/\map)=2$. There exists a regular local  subring $B$ of $R/\map$, containing a coefficient field (in the mixed characteristic case, containing a complete discrete valuation ring  whose maximal ideal is generated by the element $p^{1/2}$ of $R/\map$ and with the same residue field of $R$),  wherein $\mathbf{x}$ forms a minimal basis for its maximal ideal and $R/\map$ is module finite over $B$. Therefore, $[\Frac(R/\map):\Frac(B)]=2$, in the light of \cite[(6.5) Corollary]{HermannIkedaEquimultiplicity}.   Let $S$ be the integral closure of $B$ in $\Frac(R/\map)$. Since $B$ is a complete local ring so by \cite[Theorem 4.3.4]{HunekeSwansonIntegral} $S$ is finitely generated over $B$. Moreover, by virtue of, \cite{HochsterMcLaughlinSplitting} the inclusion $B\rightarrow S$ splits\footnote{Hence,  the Monomial Conjecture for the system or parameters $\mathbf{x}$ of $R$ follows from \cite{HochsterMcLaughlinSplitting}.}. Consequently, we may assume that $S=B\dsum I$ for some finitely generated $B$-module $I$. Since we have quadratic extension of fraction fields so we can deduce that $I$ has rank $1$. In particular, we may presume that, $I$ is an ideal of $B$ (since $S$ is a domain so $I$ is torsion-free). Using the $S_2$-property of $R$ for each prime ideal $\mathfrak{p}$ of $B$ we conclude that $I_\mathfrak{p}\cong B_\mathfrak{p}$ provided $\dep(B_\map)=\dime(B_\map)\le 1$ and $\dep(I_\map)\ge 2$ provided $\dep(B_\map)=\dime(B_\map)\ge 2$. Consequently, $I$ is a reflexive ideal of $B$ in the light of \cite[Proposition 1.4.1]{BrunsHerzogCohenMacaulay}. Now the result follows from the fact that reflexive ideals of unique factorization domains are principal.
      	  	   
      	  	      Now we are going to prove the above claim, i.e., $\text{embdim}(R)-\dime(R)\le 2$. Let us use a presentation $p^{1/2},X_2,\ldots,X_d$ for the sequence $\mathbf{x}$ in $R$ where $R=(V,p^{1/2})[[X_2,\ldots,X_d,Z_1,\ldots,Z_u]]/I$\footnote{We choose a presentation of $A$ in the mixed characteristic, because the Monomial Conjecture is open in the mixed characteristic. But the proof is characteristic free and someone can use presentation of $A$ in  equicharacteristic zero. Furthermore, we write $p^{1/2}$ instead of $p$, because we applied our  square root technique, developed in the proof of Corollary \ref{WeCanAssumeMinimalGenerator}, to have the extra assumption that $x_1,\ldots,x_d$ is a part of a minimal generating set of $\mam$.} is a homomorphic image of the regular local ring $(A,\man)=(V,p^{1/2})[[X_2,\ldots,X_d,Z_1,\ldots,Z_u]]$. So, $\man^2\subseteq (p^{1/2},Y_2,\ldots,Y_d)+I$. Set,  $I=(f_1,\ldots,f_l)$ where $l=\mu(I)$. We denote by $f^X_i$ the sum of those monomials of $f_i$ whose power of $X_i$ is non-zero for some $2\le i\le d$. Subsequently,   we set,  $f^{p^{1/2}}_i$, to be the sum of those monomials of $f_i-f^X_i$ whose coefficients are multiple of $p^{1/2}$. It follows that, $f_i^Z:=f_i-f^X_i-f_i^{p^{1/2}}\in V[[Z_1,\ldots,Z_u]]$, and that the coefficients of the monomials of $f_i^Z$ are all invertible. Now, set  $\ce{^{\ge 3}f_i^{\ Z}}$ (resp., $\ce{^{\le  2}f_i^{\ Z}}$) to be the sum of  the monomials of $f_i^Z$ of total degree greater than or equal to $3$ (resp., less than or equal to $2$). Since, $I\subseteq \man^2$ , so  it turns out that,  $\ce{^{\le 2}f^{\ Z}_i}$, is an $V$-linear combination of the elements of the form $\{Z_iZ_i:1\le i,j\le u\}$ with invertible coefficients in $V$. Otherwise, $\ce{^{\le 2}f^{\ Z}_i}$ and thence $f_i$, would have  a summand of the form $kZ_j^\alpha$ where, $\alpha\in \{0,1\}$, $1\le j\le u$ and $k\in V\backslash p^{1/2}V$. But this contradicts with $f_i\in \man^2$. In particular, $\ce{^{\le 2}f^{\ Z}_i}\in (Z_1,\ldots,Z_u)^2$. 
      	  	      
      	  	       On the other hand the fact that, $Z_iZ_j\in (p^{1/2},X_2,\ldots,X_d)+(f_1,\ldots,f_l)$ yields, $$Z_kZ_s=g_1p^{1/2}+\sum\limits_{i=2}^dg_iX_i+\sum\limits_{i=1}^lh_i\ \ce{^{\ge 3}f_i^{\ Z}}+\sum\limits_{i=1}^lh_i\ \ce{^{\le 2}f_i^{\ Z}},$$
      	  	       for each $1\le k,s\le u$ and power series  $h_i,g_i$. Thus an elementary, quite easy, computation shows that $(Z_1,\ldots,Z_u)^2\subseteq (\ce{^{\le 2}f_1^{\ Z}},\ldots,\ce{^{\le 2}f_l^{\ Z}})$. This in conjunction with the concluding assertion of the preceding paragraph yields $(Z_1,\ldots,Z_u)^2=(\ce{^{\le 2}f_1^{\ Z}},\ldots,\ce{^{\le 2}f_l^{\ Z}})$.

      	  	          Since $R$ is an almost complete intersection, so  we have, $l=\mu(I)=\hit(I)+1=\dime(A)-\dime(R)+1=u+1$.  Consequently, we get, $$u+1=l\ge\mu\Big(\sum\limits_{i=1}^{l}(\ce{^{\le 2}f_i^{\ Z}})V[[X_2,\ldots,X_d,Z_1,\ldots,Z_u]]\Big)=\mu\big((Z_1,\ldots,Z_u)^2\big)=u(u+1)/2,$$ i.e. $u\le 2$.  
      	  \end{proof}
      \end{thm}
      
      The proof of the above theorem implies also the following result which is worth stating it separately as a theorem. Recall that by \cite{RobertsAbelian} the integral closure of a regular local ring $A$ in a finite Galois extension $L$ of $K:=\Frac(A)$ with Abelian Galois group $\text{Gal}_{L/K}$ is Cohen-Macaulay provided the order of $\text{Gal}_{L/K}$ is not divisible by the characteristic of the residue field of $A$. In view of the following theorem for the quadratic extension case there is no restriction on the characteristic of the residue field, provided $A$ is a complete regular local ring.
      
      \begin{thm} \label{CohenMacaulayQuadratic}
      	Suppose that $A$ is a complete regular local ring and $L$ is a quadratic field extension of $\Frac(A)$. Then the integral closure of $A$ in $L$ is Cohen-Macaulay.
      \end{thm}
      
      \begin{rem}
        In this remark we will state some nice properties as well as several questions about the  class of rings of Theorem \ref{MonomialConjectureHolds}. So we are in the case of the  Theorem \ref{MonomialConjectureHolds} with an extra assumption that the system of parameters $\mathbf{x}$ of (\ref{StrongHypothesis}) is a part of a minimal basis of $\mathfrak{m}$  whence $\text{embdim}(R)-\dime(R)\le 2$.
          \begin{enumerate}
            \item[(i)] All of the examples of such almost complete intersections that we have computed in Macaulay2 system had $\dep\ge d-2$. So we guess  $\dep(R)\ge d-2$  whence $\omega_R$ is a maximal Cohen-Macaulay module by Lemma \ref{CanonicalModule} (ii). This observation leads us to the following question. 
              \begin{ques}\label{AlmostCompleteIntersectionQuestion}
                Suppose that $U$ is an almost complete intersection. Then do we have,
                  \begin{equation} 
                    \label{AlmostCompleteIntersectionGuess}
                   \dep(U)\ge \dime(U)-\et(U)?
                  \end{equation} 
              \end{ques}

            Question \ref{AlmostCompleteIntersectionQuestion}  is certainly true in dimension $\le 1$. But even in dimension two we need to use 	the validity of the Monomial Conjecture for the confirmation of the  inequality (\ref{AlmostCompleteIntersectionGuess}). Namely, in dimension two we only need to show that if $\dep(U)=0$ then $\et(U)\ge 2$. But if $\dep(U)=0$ then we have $\Ht_2(\mathbf{y},U)\neq 0$, wherein $\mathbf{y}$ is a system of parameters for $U$. Hence,
              $$\et(U)=\ell\big(U/(\mathbf{y})\big)-\ell\big(\Ht_1(\mathbf{y},U)\big)+\ell\big(\Ht_2(\mathbf{y},U)\big)\ge 2,$$ as the validity of the Monomial Conjecture implies that, $ \ell\big(U/(\mathbf{y})\big) -\ell\big(H_1(\mathbf{y},U)\big)\gneq 0$ by \cite[1.3. Proposition]{DuttaTheMonomial }.
              
          For another case where   inequality (\ref{AlmostCompleteIntersectionGuess})  holds, consider an almost complete intersection ideal $\maa$ of an equicharacteristic complete regular local ring $A$ with infinite residue field such that $\text{Assh}(A/\maa)=\{\map/\maa\}$. Set, $U:=A/\maa$. Suppose, furthermore, that $\dime(U)=3$ and $S:=A/\map$ is an integrally closed domain\footnote{In each prime ideal $\maq$ of $A$ we can find an almost complete intersection ideal $\maa\subseteq \maq$ such that $\text{Assh}(A/\maa)=\{\maq/\maa\}$ (see, the proof of \cite[1.2 Proposition]{DuttaTheMonomial}). So, indeed, there exists concrete examples of such almost complete intersections.}. If the inequality (\ref{AlmostCompleteIntersectionGuess}) fails for $U$ then we must have $\dep(U)=0$ and $\et(U)\le 2$.  It is easily seen that $\omega_{U}=\omega_S$. On the other hand, Lemma \ref{CanonicalModule}(ii) implies that $\dep(\omega_U)=2$ and thence $\omega_U=\omega_{S}$ is not a maximal Cohen-Macaulay module. In particular, $S$ is not a Cohen-Macaulay ring.  By \cite[(4.15) Remark.]{HermannIkedaEquimultiplicity} there exists a system of parameters $\mathbf{x}$ of $U$ such that $\et(\mathbf{x},U)=2$.  Since we are in the equicharacteristic case so there exists a complete regular local subring $B$ of $S$ with the same residue field as $S$ such that $\mathbf{x}$ forms a minimal basis of the  maximal ideal of $B$. Consequently, in view of \cite[(1.10.1.)]{HermannIkedaEquimultiplicity},  $\et(\mathbf{x},S)\le 2$. If $\et(\mathbf{x},S)=1$ then \cite[(6.8) Theorem]{HermannIkedaEquimultiplicity} implies that $S$ is regular, a contradiction. On the other hand if $\et(\mathbf{x},S)=2$, then \cite[(6.5) Corollary]{HermannIkedaEquimultiplicity} implies that $[\Frac(S):\Frac(B)]=2$. Now $S$ is Cohen-Macaulay by Theorem \ref{CohenMacaulayQuadratic}, again, a contradiction.
          
          Note that the inequality  (\ref{AlmostCompleteIntersectionGuess}) does not hold if we relax the almost complete intersection condition. More precisely, taking into account the Abhyankar's local domains \cite[(3)]{AbhyankarLocal}, we can construct local domains with multiplicity two and depth $1$ of arbitrary dimension.
              
              \item[(ii)] In view of \cite[1.3. Proposition]{DuttaTheMonomial } the Monomial Conjecture for the system of parameters $\mathbf{x}:=p^{1/2},X_2,\ldots,X_d$ of $R$ is equivalent to the assertion that, 
                $$ \ell\big(H_1(\mathbf{x},R)\big) \lneq \ell\big(R/(\mathbf{x})\big)=\ell\big(K[[Z_1,Z_2]]/(Z_1^2,Z_2^2,Z_1Z_2)\big)=3.$$
                
                On the other hand in the light of Lemma \ref{CanonicalModule}(i) the Monomial Conjecture holds for $\mathbf{x}$ if and only if $\Ht_1(\mathbf{x},R)$ is a quotient of $(Z_1,Z_2)/(Z_1^2,Z_2^2)$ by a non-zero submodule. It is easily seen that this is equivalent to the assertion that $\mathfrak{m}\Ht_1(\mathbf{x},R)=0$. Accordingly, in the light of \cite[Theorem 4.4.]{HassanzadehNaelitonResidual}, there exists a  complex, denoted by, $_0\mathcal{Z}^+_{\bullet}(\mathbf{x},Z_1)$, such that $\text{H}_0\big(_0\mathcal{Z}^+_{\bullet}(\mathbf{x},Z_1)\big)=R/\mam$ and $\Ht_1(_0\mathcal{Z}^+_{\bullet}(\mathbf{x},Z_1))=0$  if and only if the Monomial Conjecture holds for $R$. Consequently, using Theorem \ref{MonomialConjectureHolds}, we have $\Ht_1(\mathbf{x},R)$ is a non-zero at most two dimensional vector  space and $\Ht_1\big(_0\mathcal{Z}^+_{\bullet}(\mathbf{x},Z_1)\big)=0$.  Here it is noteworthy to mention that, we strongly, guess more, i.e. $\mathcal{Z}^+_{\bullet}(\mathbf{x},Z_1)$ is an acyclic finite complex consisting of Koszul cycles of $x_1,\ldots,x_d,Z_1$ which resolves the residue field of $R$ (This is equivalent to say that $\mam\Ht_i(\mathbf{x},R)=0$  for each $i\ge 1$). 
             \item [(iii)] We claim that the parameter ideal $\maa=(p^{1/2},x_2,\ldots,x_d)$ of $R$ is a Lech ideal in the sense that $\maa/\maa^2$ 	is a free $R/\maa$-module. Before proving, we wish to clear why this claim sounds interesting to us. This fact in conjunction with the exact sequence in the paragraph prior to \cite[Proposition 2.6]{CorsoHunekeIntegral} shows that $\Ht_1(\mathbf{x},R)\cong \delta_1(\maa)$ wherein, $$\delta_1(\maa)=\big(\text{Z}_1(\mathbf{x},R)\ins (\maa R^d)\big)/\text{B}_1(\mathbf{x},R)\cong \ker\ \big(\text{sym}_2(\maa)\rightarrow \maa^2\big),$$  by \cite[Corollary 2.6.]{HerzogSimisKoszul}. So \cite[XV., Proposition 12]{AndreHomologie} implies that $\Ht_1(\mathbf{x},R)$ is the second Andr\'e-Quillen homology $D_2(R,R/\maa,R/\maa)$. Hence, in view of the preceding part, it would be highly desirable to give a proof of the validity of the Monomial Conjecture for the sequence $\mathbf{x}$ of $R$ with  aid of the Andr\'e-Quillen homology avoiding, \cite{HochsterMcLaughlinSplitting}. The fact that $\maa$ is  Lech ideal follows from the following two steps.
             
             \textbf{Step 1:} Bearing in mind the notations of the proof of Theorem \ref{MonomialConjectureHolds}, this step is devoted to prove that, without loss of generality, we can assume that $f_1^Z=Z_1^2,f_2^Z=Z_2^2\ \text{and} f_3^Z=Z_1Z_2$. Recall that $R$ is a homomorphic image of  $(A,\man)=(V,p^{1/2}V)[[X_2,\ldots,X_d,Z_1,Z_2]]$, by the ideal $(f_1,f_2,f_3)$. As we have seen in the proof of Theorem \ref{MonomialConjectureHolds}, $(f^Z_1,\ldots,f^Z_3)\subseteq (Z_1,Z_2)^2$. On the other hand, by our hypothesis, $Z_iZ_j=\sum\limits_{k=2}^dX_kg_k+p^{1/2}h+\sum\limits_{k=1}^3f_k^Zl_k$, for some,  $$g_2,\ldots,g_d,h,l_1,l_2,l_3\in A.$$ We use a similar method as in the case of $f_i$ and write  $h=h^X+h^{Z,p^{1/2}}$ and $l_s=l^X_s+l^{Z,p^{1/2}}_s$ for each $1\le s\le 3$, so that $h^{Z,p^{1/2}},l^{Z,p^{1/2}}_1,\ldots,l^{Z,p^{1/2}}_3\in V[[Z_1,Z_2]].$
             Therefore,  $$Z_iZ_j-p^{1/2}h^{Z,p^{1/2}}-\sum\limits_{k=1}^3l^{Z,p^{1/2}}_kf_k^Z=\sum\limits_{k=2}^dX_kg_k+p^{1/2}h^X+
             \sum\limits_{k=1}^3f_k^Zl^X_k,$$ which implies that, $$Z_iZ_j=p^{1/2}h^{Z,p^{1/2}}+\sum\limits_{k=1}^3l^{Z,p^{1/2}}_kf_k^Z\in (f^Z_1,f^Z_2,f^Z_3,p^{1/2}).$$ Set, $B:=A/p^{1/2}=k[[X_2,\ldots,X_d,Z_1,Z_u]]$ where $k$ is the residue field of $V$. By the above argument we have $(f^Z_1,f^Z_2,f^Z_3)=(Z_1,Z_2)^2$ in $B$. Thus there exists a square matrix $H=[h_{st}]$ such that $H[f_1^Z,f_2^Z,f_3^Z]^\tau=[Z_1^2,Z_2^2,Z_1Z_2]^\tau$, in $B$. Let us to denote the maximal ideal of $B$ by $\man_B$. In particular,  $\overline{H}:=[\overline{h_{st}}]$, is the change of basis matrix with respect to the two basis $\{\overline{f_1^Z},\overline{f_2^Z},\overline{f_3^Z}\}$ and $\{\overline{Z_1^2},\overline{Z_2^2},\overline{Z_1Z_2}\}$ of the $K$-vector space $(Z_1,Z_2)^2/\man_B(Z_1,Z_2)^2$. Whence, $\det(\overline{H})\neq 0$, in $K$. Let $Y=[y_{ij}]$ be a square matrix which lifts the matrix $H$. Then $\det(Y)\notin \man$. So $Y$ is invertible. Hence, if we put, $Y[f_1,f_2,f_3]^\tau=[f'_1,f'_2,f'_3]^\tau$, then $I=(f'_1,f'_2,f'_3)$.  Moreover for each $1\le i\le 3$ we have the identity,
             \begin{equation}
             \label{FirstEqualityOfFPrime}
             f'_i=\sum_{j=1}^3y_{ij}f^X_j+\sum_{j=1}^3y_{ij}f^{p^{1/2}}_j+\sum_{j=1}^3y_{ij}f^Z_j.
             \end{equation} 
             
             On the other hand, since $H[f^Z_1,f^Z_2,f^Z_3]^\tau=[Z_1^2,Z_2^2,Z_1Z_2]$ in $B$, so we get                    \begin{equation}
             \label{SecondtEqualityOfFPrime}
             \sum\limits_{j=1}^3 y_{1j}f_j^Z\equiv Z_1^2,\text{\ } \sum\limits_{j=1}^3 y_{2j}f_j^Z\equiv Z_2^2 \text{\ \ and\ } \sum\limits_{j=1}^3 y_{3j}f_j^Z\equiv Z_1Z_2 \text{,\ \ \ mod\ } p^{1/2}A.
             \end{equation}
             So our claim follows from (\ref{FirstEqualityOfFPrime}) and (\ref{SecondtEqualityOfFPrime}).
             
             \textbf{Step 2:} In this step we assume that $\sum\limits_{i=1}^dr_ix_i=0\in R$ and we prove that, $$r_i\in \maa=(p^{1/2},x_2,\ldots,x_d),$$ for each $1\le i\le d$ ($x_1=p^{1/2}$)\footnote{This, immediately, shows that $\maa$ is a Lech ideal.}.  If so, then there are $s_1,\ldots,s_d\in (A,\man)$ such that,
              \begin{equation}
              \label{Identity}
             \sum\limits_{i=1}^d(r_i+s_i)X_i=Z_1^2h_1+Z_2^2h_2+Z_1Z_2h_3,
             \end{equation}
              because $f_1^Z=Z_1^2,f_2^Z=Z_2^2,\ \text{and} f_3^Z=Z_1Z_2$. Clearly, we must have, $h_i\in \man$, for each $1\le i\le 3$. Hence, $s_i\in \man^2$, for every $1\le i\le d$, because $f_j^X+f_j^{p^{1/2}}\in \man^2$ $(1\le j\le 3)$.  On the other hand, (\ref{Identity}) shows that, $ \sum\limits_{i=1}^d(r_i+s_i)X_i+(Z_1^2,Z_2^2,Z_1Z_2)=0$ in the Cohen-Macaulay ring $A/(Z_1^2,Z_2^2,Z_1Z_2)$ whence $r_i+s_i\in (X_1,\ldots,X_d,Z_1^2,Z_2^2,Z_1Z_2),$ for each $1\le i\le d$. Consequently, $$r_i\in (X_1,\ldots,X_d,Z_1^2,Z_2^2,Z_1Z_2)=(X_1,\ldots,X_d,f_1,f_2,f_3),$$ for each $1\le i\le d$,  because $s_i\in \man^2$.
             
          \end{enumerate}
      \end{rem}

      \begin{exam} We give an example to show that, in general, the inclusion $\mathfrak{m}^2\subseteq (\mathbf{x})$ in conjunction with the non-Cohen-Macaulayness does not imply that $\et(\mathbf{x},R)=2$ (without assuming that $R$ 	is an almost complete intersection). Set,
      	  $$ R=\mathbf{Q}[X_1,\ldots,X_4,Z_1,Z_2,Z_3] ,$$ wherein $X_1,X_2,X_3$ have degree $1$ and $X_4,Z_1,Z_2,Z_3$ have degree $2$. Let,
      	    $$I:=(Z_1^2+X_3^2Z_2+X_4Z_2,Z_2^2-X_1X_2Z_3+X_4Z_3,Z_3^2,Z_1Z_2,Z_1Z_3,Z_2Z_3).$$
      	  Then, in $S=R/I$ the image of the sequence $\mathbf{x}:=X_1,X_2,X_3,X_4$ forms a system of parameters satisfying (\ref{StrongHypothesis}) while $\et(\mathbf{x},R)=3$. Note that $\dep(R)=3$  and $R$ is not Cohen-Macaulay.
      \end{exam}
      
      \section{A question}
      We end the paper with the following question.
       
       \begin{ques}
       	  In prime characteristic, $\Ht^d_\mam(R)=\lim\limits_{\underset{n\in \mn}{\longrightarrow}}R/(\mathbf{x}^n),$ is a cyclic module over the Frobenius Skew Polynomial ring $R[x;f]$, which is  generated by $[1/(\mathbf{x})]$. This shows, in particular, that, $[1/(\mathbf{x})]\neq 0$, by the Grothendieck's non-vanishing theorem, i.e. the Monomial Conjecture holds. In characteristic zero also the Frobenius-like endomorphism exists in a faithfully flat extension with aid of  the Schoutens's Lefschetz extensions (see, \cite{AschenbrennerSchoutensLefschetz}). Hence probably the same strategy works in the equal characteristic zero. Here it seems to be natural to ask whether in the mixed characteristic there exists an $R$-algebra $S$ such that $\Ht^d_\mam(R)$ is a cyclic $S$-module generated by $[1/(\mathbf{x})]$? It is worth to point out that if this question has an affirmative answer then $S$ is necessarily a non-commutative ring. To be more precise suppose that there exists such an $R$-algebra $S$, but to the contrary $S$ is commutative. Then since $\Ht^d_\mam(R)$ 	is a cyclic $S$-module by our assumption so $\Ht^d_\mam(R)$ carries a commutative $R$-algebra structure. In particular there exists a surjective  map $\Ht^d_\mam(R)\tenr \Ht^d_\mam(R)\rightarrow \Ht^d_\mam(R)$ induced by the multiplication. But on the other hand $\Ht^d_\mam(R)\tenr \Ht^d_\mam(R)\cong \Ht^d_\mam\big(\Ht^d_\mam(R)\big)=0$ which contradicts with the Grothendieck's non-vanishing theorem (here, we assume that $d>0$).
       	  
    \end{ques}

       \section*{Acknowledgement}
       We would like to express our deepest gratitude to Massoud Tousi,  Anurag K. Singh, Linquan Ma, S. H. Hassanzadeh as well as Kazuma Shimomoto for their helpful comments and valuable suggestions. We also would like to thank the University of Utah for its hospitality during our visit in the academic year 2015-2016 as well as for its support for this research.

     \small \textsc{Ehsan Tavnafar, Department of Mathematics, Shahid Beheshti University, G.C., Tehran, Iran.}  \\
     E-mail address: \href{mailto:tavanfar@ipm.ir}{tavanfar@ipm.ir}

    \end{document}